\documentclass[acmsmall,screen]{acmart}
\pdfoutput=1

\def\BibTeX{{\rm B\kern-.05em{\sc i\kern-.025em b}\kern-.08emT\kern-.1667em\lower.7ex\hbox{E}\kern-.125emX}}

\setcopyright{acmcopyright}
\acmJournal{TOMS}
\acmYear{20xx}\acmVolume{0}\acmNumber{0}\acmArticle{0}\acmMonth{3}
\acmDOI{xx.xxxx/xxxxxx}

\usepackage{amsmath, amsthm, amssymb, amsfonts, mathtools, xfrac, mathrsfs}
\usepackage{subcaption}
\usepackage[algo2e,ruled]{algorithm2e}
\usepackage{pgfplots}
\pgfplotsset{compat=newest}
\usetikzlibrary{intersections}
\usepgfplotslibrary{fillbetween}

\usepackage[capitalize]{cleveref}
\crefname{figure}{Fig.}{Figs.}
\Crefname{figure}{Figure}{Figures}
\crefname{table}{Tab.}{Tabs.}
\Crefname{table}{Table}{Tables}
\crefname{equation}{Eq.}{Eqs.}
\Crefname{equation}{Equation}{Equations}
\crefname{section}{Sec.}{Secs.}
\Crefname{section}{Section}{Sections}
\crefname{algorithm}{Alg.}{Algs.}
\Crefname{algorithm}{Algorithm}{Algorithms}
\crefname{theorem}{Thm.}{Thms.}
\Crefname{theorem}{Theorem}{Theorems}
\crefname{corollary}{Cor.}{Cors.}
\Crefname{corollary}{Corollary}{Corollaries}

\definecolor{cerulean}{rgb}{0.0, 0.48, 0.65}
\definecolor{carrot}{rgb}{0.93, 0.57, 0.13}
\definecolor{applegreen}{rgb}{0.0, 0.5, 0.0}
\definecolor{amethyst}{rgb}{0.6, 0.4, 0.8}
\definecolor{airforceblue}{rgb}{0.36, 0.54, 0.66}
\definecolor{bittersweet}{rgb}{1.0, 0.44, 0.37}
\definecolor{citrine}{rgb}{0.89, 0.82, 0.04}
\definecolor{columbiablue}{rgb}{0.61, 0.87, 1.0}
\definecolor{flame}{rgb}{0.89, 0.35, 0.13}
\definecolor{internationalorange}{rgb}{1.0, 0.31, 0.0}

\definecolor{python0}{RGB}{31,119,180}
\definecolor{python1}{RGB}{255,127,14}
\definecolor{python2}{RGB}{44,160,44}
\definecolor{python3}{RGB}{214,39,40}
\definecolor{python4}{RGB}{148,103,189}
\definecolor{python5}{RGB}{140,86,75}
\definecolor{python6}{RGB}{227,119,194}
\definecolor{python7}{RGB}{127,127,127}
\definecolor{python8}{RGB}{188,189,34}
\definecolor{python9}{RGB}{23,190,207}

\newcommand{\Z}{\mathbb{Z}}

\newcommand{\Am}[1]{A^{(#1)}}

\newcommand{\Vmat}[1]{V_{(#1)}}

\newcommand{\probdim}{\ensuremath{N}}
\newcommand{\neig}{\ensuremath{n_e}}

\newcommand{\nsft}{\ensuremath{n_s}}

\DeclareMathOperator*{\argmin}{arg\,min}
\DeclarePairedDelimiter{\ceil}{\lceil}{\rceil}

\begin{document}

\title[Shift Selection for Parallel Shift-Invert Spectrum Slicing]{A Shift Selection Strategy for Parallel Shift-Invert Spectrum Slicing in Symmetric Self-Consistent Eigenvalue Computation}

\author{David B. Williams--Young}
\email{dbwy@lbl.gov}
\affiliation{%
  \institution{Lawrence Berkeley National Laboratory}
  \city{Berkeley}
  \state{California}
  \postcode{94720}
}

\author{Paul G. Beckman}
\email{pbeckman@uchicago.edu}
\affiliation{%
  \institution{University of Chicago}
  \city{Chicago}
  \state{Illinois}
  \postcode{60637}
}

\author{Chao Yang}
\email{cyang@lbl.gov}
\affiliation{%
  \institution{Lawrence Berkeley National Laboratory}
  \city{Berkeley}
  \state{California}
  \postcode{94720}
}

\begin{abstract}
The central importance of large scale eigenvalue problems in scientific computation
necessitates the development of massively parallel algorithms for their solution.  Recent
advances in dense numerical linear algebra have enabled the routine treatment of eigenvalue
problems with dimensions on the order of hundreds of thousands on the world's largest 
supercomputers. In cases where dense treatments are not feasible, Krylov subspace methods
offer an attractive alternative due to the fact that they do not require storage of the
problem matrices. However, demonstration of scalability of either of these classes of
eigenvalue algorithms on computing architectures capable of expressing massive 
parallelism is non-trivial due to communication requirements and serial bottlenecks, 
respectively. In this work, we introduce the SISLICE method: a parallel shift-invert
algorithm for the solution of the symmetric self-consistent field (SCF) eigenvalue problem. The
SISLICE method drastically reduces the communication requirement of current parallel
shift-invert eigenvalue algorithms through various shift selection and migration 
techniques based on density of states estimation and k-means clustering, respectively.
This work demonstrates the robustness and parallel performance of the SISLICE method on a
representative set of SCF eigenvalue problems and outlines research directions
which will be explored in future work.
\end{abstract}

\ccsdesc[500]{Mathematics of computing~Mathematical software performance}
\ccsdesc[500]{Applied computing~Mathematics and statistics}
\keywords{Eigenvalues, Parallel Eigenvalue Algorithms, Self-Consistent Field, Shift-Invert Spectrum Slicing}

\maketitle

\section{Introduction}
\label{sec:intro}

Large-scale symmetric eigenvalue problems arise in many types of scientific
computation~\cite{scidac}. In the case of electronic structure calculations based on
the Hartree--Fock approximation or Kohn-Sham density functional theory, a large symmetric nonlinear
eigenvalue problem must be iteratively solved through what is known
as the self-consistent field (SCF) procedure. Typical methods to solve the so-called
SCF eigenvalue problem require the partial diagonalization of a sequence
of matrix pencils where each pencil of the sequence is generated using
a subset of the eigenvectors of the previous pencil. The
SCF problem is considered solved when convergence of the sequence is achieved,
i.e. the change in the matrix pencil (or equivalently, the
desired eigenvectors of said pencil) between two subsequent iterations of
the SCF procedure falls below a specified threshold. 
One of the hallmarks of the SCF procedure is that the desired eigenpairs need
not be computed to full accuracy before convergence is reached. However, as the
sequence progresses, an increasing level of accuracy in the desired eigenpairs
is needed to ensure convergence to the proper SCF solution. Although the
spectrum of the matrix pencils may change quite a bit in the first few iterations of the SCF
procedure, this change becomes progressively smaller as convergence is
reached.  This feature makes the design and implementation of algorithms
for solving the SCF eigenvalue problem somewhat different from traditional
algorithms for solving an eigenvalue problem of a fixed matrix.
The repeated diagonalization
required by the SCF procedure is often the computational bottleneck in large
scale electronic structure calculations \cite{banerjee18two,laurent99electronic,shepard1993elimination}, especially in cases where a large
number of computational resources are available. As such, methods must be
developed to efficiently solve this class of nonlinear eigenvalue problem
on modern, massively-parallel computing architectures.

In this work, we consider the partial diagonalization of $\neig$ eigenpairs of
a converging sequence of symmetric matrix pencils, $(\Am{i}, B)$, of dimension \probdim,
\begin{equation}
\Am{i} X^{(i)} = B X^{(i)} \Lambda^{(i)}, \label{eq:seq_gep}
\end{equation}
where $i \in \Z^+$ is a sequence index, $\Am{i}\in\mathbb R^{\probdim \times \probdim}$ is
symmetric and $B\in\mathbb R^{\probdim \times \probdim}$ is symmetric positive definite (SPD).
$X^{(i)} \in \mathbb R^{\probdim\times\neig}$ and $\Lambda^{(i)} \in \mathbb R^{\neig\times\neig}$
are the eigenvectors and the diagonal matrix of eigenvalues corresponding to
the desired eigenpairs of $(\Am{i}, B)$, respectively.
We denote the eigenvalues as $\Lambda^{(i)}_{pq} = \delta_{pq} \lambda^{(i)}_p$ and
will refer to the increment of $i$ as an \emph{SCF iteration}.
Further, we will make the following assumptions about the sequence of matrix pencils:
\begin{itemize}
  \item We assume that $\Am{i+1}$ depends in some (possibly non-linear) way on $(X^{(i)}, \Lambda^{(i)})$.
  \item As the SCF iterations progress, we assume $\Am{i}$ converges toward a
  matrix $A$, but are not concerned with how this convergence is achieved other
  than the requirement that the convergence is not chaotic and the desired
  eigenpairs of $(\Am{i},B)$ must be computed to progressively higher accuracy as
  this convergence occurs.
  \item We assume that the desired eigenpairs of each matrix pencil in the
  sequence are contiguous within the desired spectral region bounded by
  $\lambda^{(i)}_{\min}$ and $\lambda^{(i)}_{\max}$. Remark that these 
  bounds need not correspond to the spectral bounds of $(\Am{i},B)$.
\end{itemize}

In cases where $\neig$ is relatively small ($O(<1,000)$) compared to $N$, or
when $(\Am{i},B)$ is sparse or structured, iterative algorithms such as the
implicitly restarted Lanczos algorithm \cite{lehoucq1998arpack}, the
Jacobi-Davidson algorithm \cite{sleijpen2000jacobi,mccomb07}, and the locally
optimal block preconditioned conjugate gradient (LOBPCG) algorithm
\cite{knyazev01toward} are often very effective.  
In cases where $\neig$ is a considerable
fraction of $\probdim$ or when $\neig$ is larger than thousands or tens of
thousands, iterative algorithms become less efficient partly due to the need to
solve a projected dense eigenvalue problem as a part of the Rayleigh-Ritz
procedure via some dense eigensolver.

Dense eigensolvers, such as those available in the LAPACK \cite{LAPACK},
ScaLAPACK \cite{ScaLAPACK} and ELPA \cite{marek2014elpa} libraries, are also
often used to perform a full diagonalization of each $(\Am{i},B)$. Recent
advances in dense numerical linear algebra have made it possible to perform
full diagonalizations for matrices of dimension $O(10,000) - O(100,000)$ in a
few wall clock minutes using thousands to tens of thousands computational
cores.  However, making further improvements when even more computational
resources (e.g. GPU accelerators) become available appears to be difficult due
to the communication requirement of existing parallel algorithms.

In this report, we present the SISLICE method: a parallel symmetric eigensolver
based on shift-invert spectrum slicing for the solution of the SCF eigenvalue
problem described in \cref{eq:seq_gep}. In spectrum slicing methods, the
eigenspectrum of the problem of interest is divided into several subintervals
(spectral slices) such that the eigenvalues within each slice may be computed
simultaneously. This approach eliminates the Rayleigh-Ritz bottleneck and
increases the potential for concurrency in a parallel implementation.  The
notion of spectrum slicing is well documented in the literature for many
classes of eigenvalue problems \cite{bai2000templates,saad2011numerical}. While
the basic idea behind spectrum slicing is relatively simple, its practical
implementation on large computational resources is non-trivial.
The SISLICE method has been developed to address the following practical
issues of spectrum slicing for the SCF eigenvalue problem:
\begin{itemize}
  \item Efficiently partitioning a spectral region of interest into spectral slices
  given minimal \emph{a priori} knowledge of the eigenvalue distribution,
  \item Tracking and updating spectral slice locations using emerging knowledge
  of the dynamic (but convergent) eigenvalue distribution as it changes throughout 
  the SCF procedure, and
  \item Accelerating convergence of the desired eigenpairs using knowledge of the
  approximate eigenvectors obtained from previous SCF iterations.
\end{itemize}

Fundamental to the development of spectrum slicing methods is the choice of
method for computing the approximate eigenpairs within each spectral slice.  For
this purpose, the SISLICE method employs the shift-invert subspace iteration.
Despite its simplicity, the shift-invert subspace iteration is particularly attractive
because the convergence of the method is sufficiently fast for eigenpairs in
the spectral neighborhood of a target shift. 
From a practical perspective, it is robust and
relatively simple to implement on contemporary high-performance computing
architectures.
There are a few alternatives to the shift-invert subspace iteration which include:
\begin{itemize}

  \item The shift-invert Lanczos method (SI-Lanczos). Historically, this has
  been the method of choice for spectrum slicing methods
  \cite{grimes1994shifted,zhang2007sips,campos2012strategies,mslanczos,keceli16_448,keceli18_1806}.
  However, the dependency among the sequence of linear systems the method must solve makes it 
  difficult to achieve good parallel scalability. 

  \item Polynomial filtering based methods, such as those implemented in the
  Eigenvalue Slicing Library (EVSL) \cite{li19evsl}, apply either the subspace
  \cite{banerjee16cheby,banerjee18two} or Lanczos iteration \cite{li2016thick} to
  a bandpass matrix polynomial filter which amplifies the spectral components
  associated with eigenvalues within a particular spectral slice. 
  However, in cases when the spectral slice of interest
  is not well separated from adjacent slices or when the size of the slice is
  small, polynomial filtering methods often require the construction of extremely
  high degree polynomials to ensure convergence. In practice, this can lead to
  significant performance degradation and load imbalance in spectrum slicing based 
  on polynomial filtering.

  \item Contour integral based methods \cite{ss} such as the widely adopted
  (P)FEAST method \cite{feast1,feast2,pfeast} require the solution of a number of
  complex shifted linear systems over a quadrature discretization of the spectral
  projection operator in some spectral region of interest. 
  However, choosing contours and quadrature discretizations to balance convergence and
  computational work for arbitrary spectral domains is non-trivial
  \cite{mslanczos}.
\end{itemize}

Another advantage of using a shift-invert subspace iteration is that the method
can use the approximate eigenvectors obtained from previous SCF iterations 
as the initial guess to drastically improve the convergence rate of the subspace 
iteration for the current SCF iteration. 
We note that this feature is shared by other methods such as FEAST \cite{feast1} and polynomial
filtering methods based on the subspace iteration \cite{banerjee16cheby}. 

Over the years, several schemes have been developed for parallel shift-invert
spectrum slicing. The SIPs method of \cite{zhang2007sips} utilizes a scheme
which concurrently processes overlapping spectral subdomains through dynamically
selecting shifts by tracking the convergence of the SI-Lanczos iterations
within each subdomain. The dynamic shift selection within each subdomain is
similar to the sequential shift selection scheme presented in
\cite{grimes1994shifted}.  More recently, the SIPs method has been extended to
a dynamically scheduled two-level parallel scheme targeting the SCF eigenvalue
problem (SIESTA-SIPs \cite{keceli18_1806}) which processes non-overlapping
spectral subdomains with an improved dynamic shift selection within each domain
as developed in the SLEPSc library \cite{campos2012strategies}. While robust,
this shift selection strategy does not maximize the potential for concurrency
viz independent slices as the dynamic shift selection within a spectral
subdomain still occurs sequentially.

To address this, SISLICE adopts a strategy which allows for selection of all of
the shifts which partition a spectral region of interest at the same time.  We
utilize a spectral density estimation (also known as the density of states or
DOS) to approximate the distribution of eigenvalues within the domain of
interest as discussed in \cite{li2016thick,saad2011numerical,lin2016approximating,xi18fast}.  
In addition to allowing for the
concurrent processing of all shifts simultaneously, DOS based shift selection
also allows for a more optimal initial shift selection as we
are able to place more shifts in spectral regions which have many tightly
spaced eigenvalues that are not well separated from the rest of the spectrum
and fewer shifts in regions which have isolated clusters with small radii.

To address the eigenvalue migration of the SCF eigenvalue problem, the SISLICE
method employs an eigenvalue clustering strategy to refine shift placement.
SIESTA-SIPs also adopts a similar clustering strategy for this purpose, though
there are a number of important differences with the SISLICE scheme which will be
discussed later in \cref{sec:shift}.
Because the matrix sequence in \cref{eq:seq_gep}
can change significantly in early SCF iterations when it is far from converged,
some of the selected shifts resulting from the analysis of approximate
eigenvalues in the previous SCF iteration may not be optimal. Consequently,
some of these shifts may need to be deleted and new shifts may need to be
inserted to ensure no eigenvalue is missed and all eigenvalues within the
spectral region of interest can be computed efficiently by the shift-invert
subspace iteration.  We also discuss how this can be implemented in conjunction
with a spectral slice validation scheme in \cref{sec:shift}.

The SISLICE method is designed to minimize communication overhead and thus
improve parallel scalability at the expense of performing more local
calculations.  This strategy takes into account the recent trend in high
performance computing platforms in which floating point operations have become
cheaper due to the emergence of multicore processors and accelerators while
data movement remains costly.  In SISLICE, we compute more approximate
eigenpairs than the number of eigenvalues within a spectral slice. This
redundancy does not necessarily increase time to solution if there is an
abundance of computational resources which can accommodate such redundancy.
However, it makes the validation of eigenpairs easier and more efficient to
implement.  In particular, we will show that in SISLICE it is not necessary to
check mutual orthogonality of approximate eigenvectors obtained in different
spectral slices.  As a result, our validation scheme does not require moving
vectors across different nodes or processor groups, which is often costly.
This key feature enables SISLICE to scale to very large processor counts. 

This paper is organized as follows. \Cref{sec:theory} briefly reviews the
salient aspects of shift-invert spectrum slicing and the spectral slice
validation scheme used by the SISLICE method. \Cref{sec:shift,sec:parallel}
examine the practical issues of spectrum slicing, such as shift selection,
parallel load balance, etc.,  and how the SISLICE method aims to resolve them.
\Cref{sec:numerical} provides a series of numerical experiments which exhibit
the performance and robustness of the proposed SISLICE method, and
some additional improvements to the SISLICE method which we will
implement in the future are discussed in \cref{sec:conclusion}.

\section{Shift-Invert Spectrum Slicing}
\label{sec:theory}

\begin{figure}[t]
    \centering
    \captionsetup{justification=centering}
    \begin{tikzpicture}
\begin{axis}[ 
    axis x line=center, 
    axis line style = {draw = none},
    axis y line=none, 
    ymin=-6.2, ymax=6.2,
    xmin=-6.2, xmax=6.2,
    axis equal image,
    ticks= none,
    xticklabels={},
    width=\textwidth
]

   \draw[>=latex, <-] (-6,0) -- (1,0);
   \draw[>=latex, ->] (1.8,0) -- (6,0);
   \draw[dotted,thick]      (1.1,0)   -- (1.7,0);

   \draw[] (-4,-0.1) -- (-4,0.1);
   \draw[] (-1,-0.1) -- (-1.0,0.1);
   \draw[] (2,-0.1) -- (2,0.1);
   \draw[] (4,-0.1) -- (4,0.1);

   \node[above, text=black] at ( -4.0, 0.1 ) {$\sigma_1$    };
   \node[above, text=black] at ( -1.0, 0.1 ) {$\sigma_2$    };
   \node[above, text=black] at (  2.0, 0.1 ) {$\sigma_{n_s-1}$};
   \node[above, text=black] at (  4.0, 0.1 ) {$\sigma_{n_s}$  };

  \draw [decorate,decoration={brace,amplitude=10pt,mirror},xshift=0pt,yshift=0pt]
        (-6, -0.2) -- (-4, -0.2) node [black,midway,yshift=-0.6cm] 
        {\footnotesize Slice 1};
  \draw [decorate,decoration={brace,amplitude=10pt,mirror},xshift=0pt,yshift=0pt]
        (-4, -0.2) -- (-1, -0.2) node [black,midway,yshift=-0.6cm] 
        {\footnotesize Slice 2};
  \draw [decorate,decoration={brace,amplitude=10pt,mirror},xshift=0pt,yshift=0pt]
        (2, -0.2) -- (4, -0.2) node [black,midway,yshift=-0.6cm] 
        {\footnotesize Slice $n_s$};
  \draw [decorate,decoration={brace,amplitude=10pt,mirror},xshift=0pt,yshift=0pt]
        (4, -0.2) -- (6, -0.2) node [black,midway,yshift=-0.6cm] 
        {\footnotesize Slice $n_s+1$};
  \addplot[only marks, mark=x, color=red] table {
    -5.4 0
    -5.2 0
    -5.1 0
    -4.3 0
    -4.2 0
    -4.1 0
  };

  \addplot[only marks, mark=x, color=blue] table {
    -3.9 0 
    -3.5 0
    -3.4 0
    -3.3 0
    -2.0 0
    -1.5 0
  };
  \addplot[only marks, mark=x, color=applegreen] table {
    -0.8 0
    -0.5 0
     0.5 0
  };
  \addplot[only marks, mark=x, color=amethyst] table {
    2.1 0
    2.4 0
    3.5 0
    3.6 0
    3.7 0
  };
  \addplot[only marks, mark=x, color=cerulean] table {
    4.2 0
    4.3 0
    5.0 0
  };
  \end{axis}
\end{tikzpicture}
    \caption{Partitioning the spectrum of interest into several slices or subintervals which may be computed simultaneously.}
    \label{fig:SpecPart}
\end{figure}
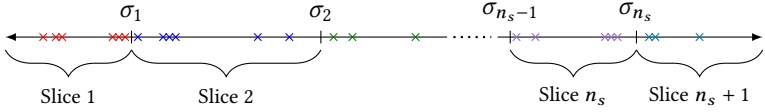

\Cref{alg:scf} outlines the major steps of the SISLICE method for solving the SCF eigenvalue problem.
At each SCF iteration, SISLICE partitions the spectral region
of interest of a matrix pencil $(\Am{i},B)$ into subintervals which may be treated independently.
These spectral subintervals will be referred to as \emph{spectral slices} in this work.
The partition is done by selecting a set of $n_s$ points
$\{\sigma_j\,\vert\,\sigma_j \in \mathbb R\}_{j=1}^{\nsft}$ as shown in \cref{fig:SpecPart}.
A particular $\sigma_j$ will be referred to as a \emph{spectral shift}.
These shift yields $n_s+1$  spectral slices.
Each spectral slice is bounded on either side by either a spectral shift or $\lambda_{\min}$ ($\lambda_{\max}$) for slice 1 ($n_s+1$), respectively.
As such, the problem of computing the eigenpairs within a particular spectral slice
amounts to computing approximate eigenpairs in the neighborhood of the shifts which form its
boundary and validating those eigenpairs against some well defined criteria.
We note here that the SISLICE method treats $n_s$ as a static quantity throughout
the SCF procedure. This constraint is adopted primarily for load balance considerations 
in a distributed parallel computing environment (see \cref{sec:parallel} for details).
In this work, we use the shift-invert subspace iteration 
to compute eigenpairs near a spectral shift.  
The approximate eigenvectors obtained from a particular
SCF iteration are used as a \emph{best guess} approximation (initial guess) for the subsequent iteration.
The validation of
eigenpairs takes into account the shifted matrix inertia as well as
the accuracy of the computed eigenpairs.
Specific details regarding the selection of
spectral shifts are given in \cref{sec:shift}. In this section, we review the salient
aspects of the shift-invert subspace iteration and slice validation scheme used by the
SISLICE method given a set of spectral shifts.

We should note that the algorithm outlined in \cref{alg:scf} can be used
to compute desired eigenpairs of a fixed matrix pencil also. When
the matrix pencil $(A,B)$ is fixed, we obviously do not need to perform the
update in Step 5, however; we may improve the efficiency of
the shift-invert subspace iteration by repartitioning the spectral region
of interest (Step 1) using previously computed eigenvalue approximations
as the reference. The application of \cref{alg:scf} to a fixed eigenvalue
problem is particularly attractive for large, sparse matrix pencils with
$N > O(100,000)$ where sparse matrix factorizations are possible, but
direct eigenvalue decomposition is impractical on currently available
computer hardware. See \cref{sec:scaling} for such an example.

\begin{algorithm2e}[t]
  \caption{Shift-invert spectrum slicing for computing $\neig$ eigenpairs of a sequence of matrix pencils $(\Am{i},B)$
  as they converge to a pencil $(A,B)$.}
  \label{alg:scf}
  \SetKwInOut{kinput}{Input}
  \SetKwInOut{koutput}{Output}
  \BlankLine
  \kinput{$(\Am{0},B)$, number of desired eigenpairs $\neig$, number of slices $K$}
  \BlankLine
  \koutput{$X \in \mathbb{R}^{\probdim\times\neig}$ and diagonal matrix $\Lambda \in \mathbb R^{\neig\times\neig}$
           which describe the desired $\neig$ eigenpairs of the converged $(A,B)$.
  }
  \BlankLine
  \BlankLine
  \BlankLine
  \For{$i = 0,1,2,\ldots$} {
    \nl Partition the spectral region of interest for $(\Am{i},B)$ into $K$ slices by selecting $n_s = K-1$ spectral shifts; \\
    \nl Choose starting guesses for eigenvectors within each slice; \\
    \nl Use shift-invert subspace iteration to obtain approximate
        eigenpairs within each slice; \\
    \nl Validate the computed eigenvalues $\rightarrow (X,\Lambda)$;\\
    \nl Compute the matrix $\Am{i+1}$ using $(X,\Lambda)$;\\
    \If{ $\Am{i+1}$ converged }{
      \nl \Return $(X,\Lambda)$;
    }
  }
\end{algorithm2e}

\subsection{The Shift--Invert Subspace Iteration}
\label{sec:sisub}

Approximate eigenpairs of $(A,B)$ corresponding to eigenvalues in the neighborhood
of a shift, $\sigma$, may be obtained through the shift-invert subspace iteration \cite{bai2000templates,saad2011numerical}.
For the purposes of this work, we assume that $\sigma$ is distinct from any eigenvalue of $(A,B)$.
One possible implementation of the shift-invert subspace iteration is given 
in \cref{alg:sisubit}.  The major cost of this implementation is in the Bunch-Kaufman (LDL$^T$)
factorization of $A-\sigma B$. When $A$ and $B$ are dense, we employ the 
LDL$^T$ factorization implemented in LAPACK or ScaLAPACK.
When $A$ and $B$ are sparse, we may use a symmetric sparse solver such as those
implemented in MUMPS~\cite{MUMPS:1,MUMPS:2},
PARDISO~\cite{schenk2000efficient,schenk2002two,schenk2006fast},
symPACK~\cite{bachan2017upc++,bachan2019upc++}, etc.  The diagonal factor $D$
produced by this factorization may be further used in the validation of
eigenpairs within a spectral slice (see \cref{sec:valid} for details).  We use
the Cholesky QR algorithm to $B$-orthonormalize the basis vector $V$ of the
subspace produced by the subspace iteration. 
This procedure requires computing the metric Grammian $G = V^T B V$, performing
a Cholesky factorization $G = RR^{T}$ and using back substitution to update $V$
via $V \leftarrow VR^{-T}$. Upon the completion of the subspace iteration, we
use the Rayleigh Ritz procedure to retrieve approximate eigenvalues and
eigenvectors of $(A,B)$ from the subspace $V$. This procedure requires solving
the projected eigenvalue problem $(V^T A V) Q = Q\Omega$ and updating $V$ by
$V\leftarrow VQ$.

Remark that the algorithm outlined in \cref{alg:sisubit} is the simplest version of
the shift-invert subspace iteration. Several modifications can be made to improve
the performance and robustness of the algorithm, especially in the context of
a convergent sequence of matrix pencils which must be treated in \cref{eq:seq_gep}.
For example, the number of columns of in $V_{(m)}$ needs only
be \emph{at least} $k$ to obtain approximations for $k$ eigenpairs. In practice,
the convergence rate of the subspace iteration in obtaining the $k$ \emph{desired}
eigenpairs may be drastically improved by choosing a trial vector space which is
several times larger than $k$ (see \cref{sec:probe_dim_results} for examples).
Further, we may exploit the fact that $(X^{(i)}, \Lambda^{(i)})$ need not be computed to full accuracy
until the SCF iterations of \cref{eq:seq_gep} are nearly converged.
The rate at which $V_{(m)}$ approaches $X$ in \cref{alg:sisubit}
largely depends on the choice of initial guess $V_{(0)}$. If the distance
between $V_{(0)}$ and $X$ (as measured in terms of subspace angle) is sufficiently small, convergence may be achieved in only a few subspace iterations.
As $\Am{i}$ converges to $A$, the change in the eigensystem between $(\Am{i},B)$ and $(\Am{i+1},B)$
becomes sufficiently small such that the distance between $X^{(i)}$ and $X^{(i+1)}$ is also small.
Thus, the subspace iteration may be seeded with $X^{(i)}$ to obtain $X^{(i+1)}$
to enable faster convergence. This assumption is typically
most valid in the last few SCF iterations, though this seeding procedure will be demonstrated
to be effective throughout the SCF procedure for the problems considered in \cref{sec:numerical}.

\begin{algorithm2e}[t]
  \caption{The Shift-Invert Subspace Iteration: $\mathrm{SISubIt}(A,B,V_{(0)},\sigma,m)$}
  \label{alg:sisubit}
  \SetKwInOut{kinput}{Input}
  \SetKwInOut{koutput}{Output}
  \BlankLine
  \kinput{Symmetric matrices $A,B \in \mathbb R^{\probdim\times\probdim}$ with $B$ being SPD,
          a target shift $\sigma \in \mathbb R$,
          the number of eigenpairs to be computed $k$,
          an initial guess of the eigenvectors $V_{(0)} \in \mathbb{R}^{\probdim\times k}$,
          and a number of subspace iterations $M$}
  \BlankLine
  \koutput{$(X,\Lambda)$ which approximates $k$ eigenpairs of $(A,B)$ in the spectral neighborhood of $\sigma$.}
  \BlankLine
  \BlankLine
  \nl $\Vmat{0} \leftarrow \mathrm{CholeskyQR}(\Vmat{0},B)$;\\
  \nl $(L,D) \leftarrow LDL^T$ factorization of $A-\sigma B$; \\
  \For{m = 1,2,...,M} {
     \nl $\Vmat{m} \leftarrow L^{-T} D^{-1}L^{-1}B\Vmat{m-1}$; \\
     \nl $\Vmat{m} \leftarrow \mathrm{CholeskyQR}(\Vmat{m},B)$; \\
  }
  \nl \Return $(X,\Lambda) \leftarrow \mathrm{RayleighRitz}(A,B,\Vmat{M})$; \\
  \BlankLine
\end{algorithm2e}

Once a Rayleigh-Ritz procedure has been performed, the convergence of the subspace iteration 
can be assessed by computing a set of residuals, $R_j = AX_j - BX_j\Lambda_j \in \mathbb{R}^{N\times k}$, and evaluating the 2-norm of each column of $R_j$. 
The tuple $(\sigma_j,\Lambda_j,X_j,R_j$) will be
referred to as the $j$-th spectral probe throughout the remainder of this work and will
be occasionally denoted SP($\sigma_j$). We note for clarity that one need not consider both
$V_j$ and $X_j$ simultaneously due to the fact that they admit identical linear spans ($X_j$
is simply a rotation of $V_j$). As $X_j$ contains more useful information related to the
eigensystem of $(A,B)$ than $V_j$, $V_j$ is typically discarded in favor of $X_j$ for eigenvalue
calculations.

\subsection{Validation of Spectral Slices}
\label{sec:valid}

In the SISLICE method, the approximate eigenpairs associated with a particular slice $(\sigma_j,\sigma_{j+1})$
are obtained by analyzing the Ritz pairs that are computed from the spectral probes defined by the spectral shifts $\sigma_j$ and $\sigma_{j+1}$, which we denote by SP($\sigma_j$) and SP($\sigma_{j+1}$), respectively.
  The Ritz values obtained from SP($\sigma_j$) can potentially overlap with those obtained from SP($\sigma_{j+1}$). It is also possible that
$\sigma_j$ and $\sigma_{j+1}$ are too far apart such that a number of desired eigenvalues are not captured by either SP($\sigma_j$) or SP($\sigma_{j+1}$). Thus,
a key aspect in the development of a robust spectrum slicing method is to provide a mechanism
to select approximate eigenpairs within a spectral slice from the Ritz pairs of its associated probes
as to avoid double counting and detect any missing or spurious eigenpairs, if present. Such selected
eigenpairs will be referred to as being \emph{validated}.

\begin{figure}[t]
\centering
\begin{tikzpicture}
\begin{axis}[
    axis x line=center,
    axis y line=none,
    ymin=-4, ymax=4,
    xmin=-4, xmax=4,
    axis equal image,
    xtick={ -2, 2 },
    xticklabels={},
    every tick/.style={
        black,
        thick,
    },
    width=0.8\textwidth
]
\draw[>=latex, <->] (-4,0) -- (4,0);
\draw[dashed]       (0, -0.5) -- (0, 0.5);
\node[below,text=black] at (-2, 0 )     {$\sigma_j$};
\node[above,text=black] at (2, 0 )      {$\sigma_{j+1}$};
\addplot [only marks,
          fill opacity=0,
         ] table {
  -3.8 0.15
  -3.0 0.15
  -2.3 0.15
  3.8 -0.15
  3.0 -0.15
  2.1 -0.15
};
\addplot [only marks] table {
  -1.9 0.15
  -1.0 0.15
  -0.7 0.15
  -0.4 0.15
  -0.2 0.15
  1.9 -0.15
  1.7 -0.15
  0.4 -0.15
};
\addplot [only marks,
          mark=triangle*,
          color = red
  ] table {
  3.95 0.15
  3.1 0.15
  0.4 0.15
  -3.6 -0.15
  -1.2 -0.15
  -0.8 -0.15
  -0.4 -0.15
  -0.2 -0.15
};

\end{axis}
\end{tikzpicture}
\caption{Scheme for selection of validation candidates for the spectral slice bounded by $(\sigma_j,\sigma_{j+1})$.
The points with markers above the axis are Ritz values computed from the
SP($\sigma_j$), and those below the axis are those computed from
SP($\sigma_{j+1}$). The vertical dashed line denotes the midpoint
of the spectral slice ($\tau$). The filled black circles represent the validation candidates for
the slice, while the open circles represent Ritz values which may belong to other
spectral slices depending on the placement of $\sigma_{j-1}$ and $\sigma_{j+2}$. The
red triangles represent Ritz values which are not considered for validation for the spectral slice.
}
\label{fig:evalidate}

\end{figure}
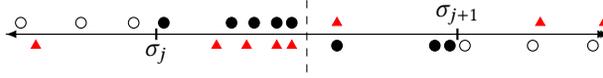

To select candidates for validation within the spectral slice
$(\sigma_j,\sigma_{j+1})$, we examine the Ritz values computed from the
spectral probes SP($\sigma_j$) and SP($\sigma_{j+1}$) that are within
$(\sigma_j,\sigma_{j+1})$. We choose a point $\tau$ between $\sigma_j$ and
$\sigma_{j+1}$, e.g., $\tau=(\sigma_j+\sigma_{j+1})/2$ and select all Ritz
values obtained from SP($\sigma_j$) that are in $(\sigma_j,\tau)$, and those
from SP($\sigma_{j+1}$) that are in $(\tau,\sigma_{j+1})$ as validation
candidates.  A graphical representation of this candidate selection process is
depicted in \cref{fig:evalidate}. For the spectral slices at both ends of the
desired spectral region of interest, validation candidates are selected as
those Ritz values that are in $[\lambda_{\min},\sigma_1)$ and
$(\sigma_{n_s+1},\lambda_{\max}]$, respectively.  This partitioning strategy
follows from the assumption that Ritz values that approximate eigenvalues
closer to a spectral shift tend to converge faster. That is, if $\theta_j$ and
$\theta_{j+1}$ are both approximations to the same eigenvalue $\lambda$ that
lies in $(\sigma_j,\tau)$, but are obtained from two different spectral probes
SP($\sigma_j$) and SP($\sigma_{j+1}$), the residual norm
associated with $\theta_j$ is likely to be smaller because $\lambda$ is closer
to $\sigma_j$ than to $\sigma_{j+1}$. This feature is demonstrated numerically
for each of the numerical experiments examined in \cref{sec:numerical}. As a
result, this partitioning strategy is responsible for removing the majority of
eigenpair duplication between adjacent probes. We denote the number of validation
candidates for the slice bounded by SP($\sigma_j$) and SP($\sigma_{j+1}$) as
$n^{cand}_j$ in the following.

The duplication of eigenvalues can also be checked
by measuring the mutual orthogonality of the corresponding eigenvectors.
However, such a scheme would require comparing Ritz vectors computed by different
spectral probes. In a parallel implementation in which different spectral
probes are mapped to different processor groups, this scheme would require
excessive data communication.
We choose to check duplication or missing eigenvalues by
simply comparing the number of validation candidates with
the exact eigenvalue count that can be obtained from
the factorization $L_jD_jL_j^T=\Am{i}-\sigma_jB$ for each spectral shift. By
making use of Sylvester's inertia theorem \cite{sylvester1852xix}, we are able to ascertain the exact
number of eigenvalues within the slice bounded by $(\sigma_j,\sigma_{j+1})$ by
taking the difference between the number of negative diagonal elements of
$D_{j+1}$ and $D_j$, respectively. We denote this exact count as $n^{exact}_j$.

Given $n^{exact}_j$ and $n^{cand}_j$ for a particular spectral slice, there are
3 limiting cases which allow us to detect missing or spurious eigenpairs within
that slice:
\begin{itemize}
\item If $n^{exact}_j = n^{cand}_j$ we view each of the candidates as a reasonable
approximation to a true eigenpair and consider it to be validated.

\item If $n^{exact}_j > n^{cand}_j$, there are true eigenpairs within
this spectral slice that are not captured by either SP($\sigma_j$) or
SP($\sigma_{j+1})$. Thus, either more shift-invert subspace iterations need to
performed or a spectral shift must be added somewhere within
$(\sigma_j,\sigma_{j+1})$ to ensure that all desired eigenpairs are accounted
for. We examine the specifics of this shift addition in \cref{sec:misev}.

\item If $n^{exact}_j < n^{cand}_j$, a number of the validation candidates are 
either spurious eigenpairs or very poor duplicate approximations to eigenpairs which
have true eigenvalues in another spectral slice. In either of these cases,
the residual norm for these unwanted eigenpairs will be larger than the desired 
ones. As a result, they may be removed by selecting the $n^{exact}_j$
eigenpairs in the slice from the bounding spectral probes with the smallest
residual norms to be the validated eigenpairs.
\end{itemize}

\section{Shift Selection and Migration}
\label{sec:shift}

In this section, we examine the selection of the set of spectral shifts $\{\sigma_j\}_{j=1}^{n_s}$ which
partition the spectral region of interest into $n_s+1$ spectral slices. 
We first utilize a Lanczos approximation for the so-called density of states (DOS) of
to select an initial set of shifts for the first SCF iteration. In subsequent SCF iterations,
we use a clustering algorithm select shifts based on approximate eigenvalues obtained 
in the previous SCF iteration.

%
%
%

\subsection{Shift selection via Lanczos-based density of states estimation}
\label{sec:dos}

Consider the eigenvalue decomposition of the pencil $(A,B)$ in \cref{eq:seq_gep}.
The exact density of states (DOS) for $(A,B)$ is given by
\begin{equation}
\phi(\omega) = \sum_{j=1}^N \delta(\omega - \lambda_j),
\end{equation}
where $\delta(\cdot)$ is the Dirac delta distribution. 
The cumulative density of states (CDOS) which returns the
number of eigenvalues of $(A,B)$ below a certain value $\omega$ is defined as
\begin{equation}
\Phi(\omega) \equiv \int_{-\infty}^{\omega} \phi(\omega') \mathrm{d}\omega'. \label{eq:cdos_general}
\end{equation}
It follows from \cref{eq:cdos_general} that $\gamma(a,b) = \Phi(b) - \Phi(a)$
provides an eigenvalue count for $(A,B)$ on $[a,b]$.
Clearly, the construction of the exact DOS and CDOS requires a full diagonalization of
$(A,B)$, which is something we are trying to avoid. 
One way to obtain an approximate DOS is to use the Lanczos algorithm as
explained in \cite{saad2011numerical,lin2016approximating,xi18fast}. For a
generalized eigenvalue problem, one can either run the Lanczos algorithm on
$L^{-1}AL^{-T}$, where $L$ is the Cholesky factor of $B$, i.e., $B=LL^T$, or
run a $B$-orthonormal Lanczos on $B^{-1}A$. The latter requires the initial
vector of the Lanczos iteration to be generated to have an appropriate random
distribution. In this work, we select the elements of the initial starting
guess using the normal distribution $\mathcal{N}(0,1)$, though we note that
more robust schemes exist for generalized eigenvalue problems \cite{xi18fast}.

Let $\{ (\theta_j, g_j)\}_{j=1}^k$ be the Ritz pairs produced from a $k$-step
Lanczos iteration associated with a fixed random starting guess. An estimated DOS may be written as
\begin{equation}
\phi(\omega) = 
  \frac{N}{\sqrt{2\pi}}\sum_{j=1}^k \frac{\zeta_j^2}{\nu_j} \exp \left( -\kappa_j(\omega)^2 \right), \qquad \zeta_j = e_1^Tg_j, 
  \quad \kappa_j(\omega) = \frac{\omega - \theta_j}{\nu_j \sqrt{2}}.
\label{eq:dos}
\end{equation}
$\nu_j$ is a length parameter which determines the width of the Gaussian. 
For the purposes of this
work, we choose $\nu_j$ so that each Gaussian factor in the sum nearly vanishes some distance $d_j$ away from $\theta_j$.
The parameter $d_j$ is chosen to be either the maximum or average of
$\theta_j-\theta_{j-1}$ and $\theta_{j+1} -\theta_j$. Some safeguard is
used to prevent $d_j$ from becoming too small when eigenvalues are tightly clustered. 
We may obtain a closed form expression for the corresponding CDOS as
\begin{equation}
\Phi(\omega) = 
  \frac{N}{2} \sum_{j=1}^k \zeta^2_j \left[ 
    \mathrm{erf}(\kappa_j(\omega)) + 1
  \right]
\label{eq:cdos}
\end{equation}
where $\mathrm{erf}(\cdot)$ is the error function. 
To get a more accurate estimate of the DOS and CDOS, multiple Lanczos runs with
different random starting guesses may need to be used. We refer readers to 
\cite{lin2016approximating,xi18fast} for more details.
The remainder of this subsection will be dedicated to how to partition a spectral
region of interest into reasonable subintervals and select shifts for the spectral probes based
on the DOS and CDOS of $(A,B)$ produced by the Lanczos algorithm. As such, we expect
that these methods are able to produce a specified number of intervals, $n_s$.
Without loss of generality, we consider an arbitrary spectral domain $[a,b]$
with approximate eigenvalue count $C(a,b) \equiv \ceil{\gamma(a,b)}$.  

In the case where the CDOS increases gradually in a nearly
continuous fashion (e.g. \cref{fig:graphene}), we may partition the spectral
region of interest into $n_s =\frac{C(a,b)}{K}$ intervals containing roughly $K$ eigenvalues by
determining the roots of the CDOS at evenly spaced intervals
\begin{equation}
  \Phi(\omega) - (a + Kj) = 0, \qquad j \in \left[ 0, n_s \right]. \label{eq:cdos_root}
\end{equation}
As $\Phi(\cdot)$ does not admit an analytic inverse, these roots may be found
either through bisection or Newton's method.  Letting $l_j$ and $u_j$ be the
$(j-1)$-st and $j$-th root of \cref{eq:cdos_root}, 
the spectral intervals produced by this method may be given as $I = \{
(l_j, u_j) \}_{j=1}^{n_s}$.  This simple strategy is guaranteed to produce exactly a
specified number of intervals but is sub-optimal when the spectral region of
interest is highly clustered with large gaps between clusters.  These clusters
may be identified by sharp peaks in the DOS (e.g. \cref{fig:silane}).  The
previously described method will often miss these isolated clusters due to the
sampling nature of the methods used to solve \cref{eq:cdos_root}.  In this
scenario, we introduce an alternative method which aims to find these spectral
clusters via identification of the local maximizers of the DOS itself. 

\begin{algorithm2e}[t]
  \caption{Obtain intervals containing eigenvalue clusters: DOSCluster($T,[a,b],n_\omega$)}
  \label{alg:dos_part}
  \SetKwInOut{kinput}{Input}
  \SetKwInOut{koutput}{Output}
  \SetKwRepeat{Do}{do}{while}
  
  \kinput {
  Ritz pairs $T = \{(\theta_j), g_j\}_{j=1}^k$ from a $k$-step Lanczos iteration,\\
  Search interval $[a,b]$ and number of discretization points $n_\omega$.
  }
  \koutput {
   Eigenvalue cluster intervals  $I$.
  }

  \BlankLine
  \nl Let $S = \{\omega_i\}_{i=1}^{n_\omega}$ with $\displaystyle \omega_i = a + \frac{b - a}{n_\omega - 1} (i-1) $; \\
  \BlankLine
  \nl Identify all local maximizers $\{ \hat\omega_j \}_{j=1}^{n_c}$ of $\phi(\cdot)$ restricted to $S$;\\
  \BlankLine
  \nl \lIf{$n_c = 1$} {
    \Return $I = \{(a,b)\}$
  }
  \BlankLine
  \nl Identify local minimizers between maximizers, $\{\mu_j\}_{j=0}^{n_c}$ with $\mu_0=a$, $\mu_{n_c}=b$ and
  \begin{equation*} 
    \mu_j = \argmin_{\omega\in (\hat\omega_j, \hat\omega_{j+1}) \subset S} \phi(\omega), \quad j \in [1,n_c);
  \end{equation*}
  \BlankLine
  \nl Let $I = \{(\mu_{j-1},\mu_j)\}_{j=1}^{n_c}$; \\
  \nl Remove intervals in $I$ which do not contain a Ritz value;\\
  \Return I;
\end{algorithm2e}

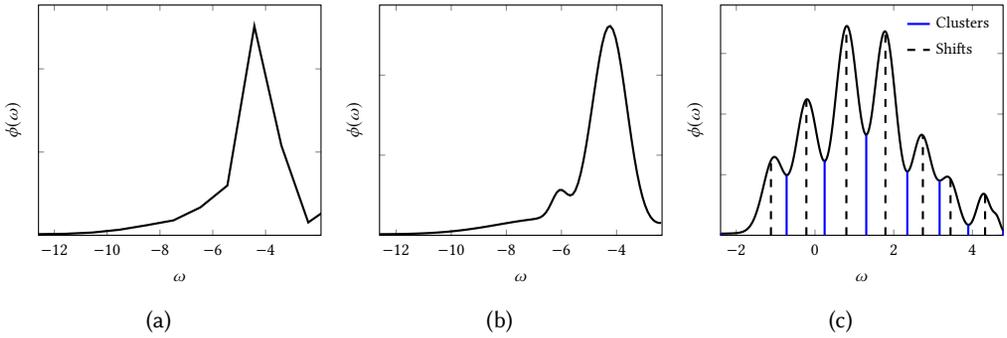
\begin{figure}[t]
  \centering
  \begin{subfigure}[b]{0.32\textwidth}
    \centering
    \begin{tikzpicture}
\pgfmathsetlengthmacro\MajorTickLength{
  \pgfkeysvalueof{/pgfplots/major tick length} * 0.5
}
\begin{axis}[
  width=1.2\textwidth,
  every axis plot/.append style={ thick},
  legend style={font=\tiny},
  legend cell align={left},
  every tick label/.append style={font=\tiny},
  major tick length=\MajorTickLength,
  xmin=-12.6, xmax=-1.9,
  ymin = 0,
  xlabel={$\omega$},
  x label style={ font=\tiny },
  xtick={ -12, -10, -8, -6, -4 },
  ylabel={$\phi(\omega)$},
  y label style={ font=\tiny },
  ylabel shift = -4 pt,
  ytick={},
  yticklabels={},
]
\addplot[mark={}] table [x=x, y=y, col sep=comma, each nth point=10] {figs/fig5Left.txt};
\end{axis}
\end{tikzpicture}
    \caption{}
    \label{fig:dos_cluster_coarse}
  \end{subfigure}
  \begin{subfigure}[b]{0.32\textwidth}
    \centering
    \begin{tikzpicture}
\pgfmathsetlengthmacro\MajorTickLength{
  \pgfkeysvalueof{/pgfplots/major tick length} * 0.5
}
\begin{axis}[
  width=1.2\textwidth,
  every axis plot/.append style={ thick},
  legend style={font=\tiny},
  legend cell align={left},
  every tick label/.append style={font=\tiny},
  major tick length=\MajorTickLength,
  xmin=-12.6, xmax=-2.35,
  ymin = 0,
  xlabel={$\omega$},
  x label style={ font=\tiny },
  xtick={ -12, -10, -8, -6, -4 },
  ylabel={$\phi(\omega)$},
  y label style={ font=\tiny },
  ylabel shift = -4 pt,
  ytick={},
  yticklabels={},
]
\addplot[mark={}] table [x=x, y=y, col sep=comma] {figs/fig5Right.txt};
\end{axis}
\end{tikzpicture}
    \caption{}
    \label{fig:dos_cluster_fine}
  \end{subfigure}
  \begin{subfigure}[b]{0.32\textwidth}
    \centering
    \begin{tikzpicture}
\pgfmathsetlengthmacro\MajorTickLength{
  \pgfkeysvalueof{/pgfplots/major tick length} * 0.5
}
\pgfplotstableread{figs/dos_spec_part_bounds.dat}{\boundstable}
\pgfplotstableread{figs/dos_spec_part_shifts.dat}{\shiftstable}

\pgfplotstablegetrowsof{figs/dos_spec_part_bounds.dat}
  \pgfmathsetmacro{\boundsrows}{\pgfplotsretval-1}
\pgfplotstablegetrowsof{figs/dos_spec_part_shifts.dat}
  \pgfmathsetmacro{\shiftsrows}{\pgfplotsretval-1}

\begin{axis}[
  width=1.2\textwidth,
  legend style={font=\tiny, draw=none, at={(1,1)}, anchor=north east, fill=none},
  legend cell align={left},
  legend image post style={xscale=0.5},
  every axis plot/.append style={ thick},
  legend style={font=\tiny},
  legend cell align={left},
  every tick label/.append style={font=\tiny},
  major tick length=\MajorTickLength,
  xmin=-2.4, xmax=4.78,
  ymin = 0,
  xlabel={$\omega$},
  x label style={ font=\tiny },
  ylabel={$\phi(\omega)$},
  y label style={ font=\tiny },
  ylabel shift = -4 pt,
  ytick={},
  yticklabels={},
]

\addlegendentry{Clusters};
\addlegendimage{line legend, no markers, blue};

\addlegendentry{Shifts};
\addlegendimage{line legend, no markers, dashed, black};

\addplot[thick,mark={},name path=dos,forget plot] table [x=x, y=y, col sep=comma] {figs/dos_spec_part_dos.dat};


\addplot[thick,ycomb, color=blue] table [x=x, y=y, col sep=comma] {figs/dos_spec_part_bounds.dat};
\addplot[thick,ycomb, color=black, dashed] table [x=x, y=y, col sep=comma] {figs/dos_spec_part_shifts.dat};

\end{axis}
\end{tikzpicture}
    \caption{}
    \label{fig:dos_cluster_shifts}
  \end{subfigure}
  \caption{Refining the resolution of the DOS (a) in the third cluster shown in
  \cref{fig:silane_dos} reveals an additional cluster near -6.0 (b). The
  refinement also allows us to set tighter bounds for refined clusters. 
  (c) depicts the interval definitions and shift placement for the upper end of \cref{fig:silane_dos}.
  } 
  \label{fig:dos_cluster}
\end{figure}

\Cref{alg:dos_part} depicts the method used to partition the DOS into intervals
containing eigenvalue clusters.  However, due to the highly non-linear nature
of of the Gaussian approximation in \cref{eq:dos} and the possibility of
spurious Ritz pairs in the $k$-step Lanczos procedure, care must be taken in
identification of the local maximizers. In particular, local maximizers must be
identified using an appropriate spectral resolution.  Evaluating DOS on a very
fine spectral grid may produce too many ``artificial" local maximizers that are
introduced by the inexact nature of the estimated DOS. Evaluating DOS on a very
coarse grid may result in missing an important local maximizer (hence a
cluster, see e.g. \cref{fig:dos_cluster_coarse,fig:dos_cluster_fine}). 
To complicate matters further, there is
no simple way of determining the ``right" discretization \emph{a priori}
to overcome these issues. To address this, we have adopted a refinement
strategy to minimize the chance of missing important spectral information
contained in the DOS approximation: 

\begin{enumerate}
  \item Obtain an initial set of cluster intervals ($I$) from \cref{alg:dos_part} on
  a coarse discretization of $[a,b]$. We have found that $n_\omega=10k$
  suffices for this purpose in most cases.
  
  \item If $\gamma(l_j,u_j)$ for some $(l_j,u_j) \in I$ is less than some specified
  tolerance, e.g. 2, merge this interval into the adjacent interval with the fewest estimated eigenvalues if it exists 
  (i.e. $(l_{j\pm1}, u_{j\pm1})$ if $u_{j-1}=l_j$ or $u_j = l_{j+1}$, respectively). If 
  there is no adjacent interval (i.e. the interval represents a cluster isolated from the rest of the spectrum),
  the interval must remain unaltered to ensure that it is captured in the shift-invert
  subspace iteration.

  \item Rather than refine the search for clusters across all of the intervals in $I$,
  it is best to limit this refined search only to intervals with large eigenvalue
  counts. Thus if $\gamma(l_j,u_j)$ is larger than some specified tolerance, e.g.
  50, we rerun \cref{alg:dos_part} with $(l_j,u_j)$ as the search interval with 
  a finer discretization. We have found that doubling the resolution within the
  spectral interval, i.e. $n_\omega' = \frac{2C(l_j,u_j)}{C(a,b)}n_\omega$, suffices
  in most cases. 

  \item If the refined search produces two or more intervals, these intervals replace
  $(l_j,u_j)$ in $I$. Steps (2) and (3) are then repeated until no new intervals are produced.

\end{enumerate}

Once the set of spectral intervals have been produced by either of the aforementioned
schemes, a shift may be selected for each interval by one of several schemes. The 
simplest scheme takes the shift to be the midpoint of the interval, e.g. 
$\sigma_j = \frac{u_j - l_j}{2}$. However, in cases when the eigenvalues are irregularly
distributed throughout the interval, we may take the expected value of $\omega$ with respect to the DOS over 
the interval to be the shift,
\begin{equation}
\sigma_j = \frac{1}{\gamma(l_j,u_j)} \int_{l_j}^{u_j} \omega\phi(\omega)\,\mathrm{d}\omega = \frac{N}{\Phi(u_j)-\Phi(l_j)}\sum_{i=1}^k \zeta_i^2 
  \left(\psi_i(u_j) - \psi_i(l_j)\right),
\end{equation}
where
\begin{equation}
\psi_i(\omega) = \frac{\theta_i}{2}\mathrm{erf}\left( \kappa_i(\omega) \right) -
\frac{\nu_i}{\sqrt{2\pi}}\exp \left( -\kappa_i(\omega)^2 \right)
\end{equation}
An example of this shift selection scheme combined with the identification of 
clusters though iterative refinement of the DOS is shown in \cref{fig:dos_cluster_shifts}.

Unlike the CDOS root finding method, \cref{alg:dos_part} is not guaranteed to
produce a specified number of intervals even after refinement. In fact, this
method produces something akin to the natural clustering of the eigenvalue
distribution. A number of safeguards have to be put in place to ensure that we
are able to produce exactly a specified number of intervals

\begin{itemize}

  \item If $\vert I \vert < n_s$, we subdivide the intervals with the largest 
  number of approximate eigenvalues.  
  This process is repeated until $\vert I \vert = n_s$.

  \item $\vert I \vert > n_s$, we merge intervals with only a few number of
  eigenvalues into adjacent intervals until $\vert I \vert = n_s$. Care must
  be taken in this procedure as merging well separated clusters (even if adjacent)
  will likely lead to poor convergence in the shift-invert subspace iterations.

\end{itemize}

\subsection{Shift Refinement and Eigenvalue Clustering}
\label{sec:kmeans}

\begin{figure}[t]
\centering
\begin{minipage}{0.7\textwidth}
\centering
\begin{tikzpicture}
\begin{axis}[ 
    axis x line=center, 
    axis line style = {draw = none},
    axis y line=none, 
    ymin=-4, ymax=4,
    xmin=-4, xmax=6,
    axis equal image,
    ticks= none,
    width=\textwidth
]

\draw[>=latex, <->] (-4,2) -- (4,2);
\draw[>=latex, <->] (-4,1) -- (4,1);
\draw[>=latex, <->] (-4,0) -- (4,0);
\draw[>=latex, <->] (-4,-1) -- (4,-1);

\node[right, text=black] at (4.5,2)  {\small DOS};
\node[right, text=black] at (4.5,1)  {\small Update 1};
\node[right, text=black] at (4.5,0)  {\small Update 2};
\node[right, text=black] at (4.5,-1) {\small Update 3};

\draw[] (-1.75, 1.9 ) -- (-1.75, 2.1);
\draw[] (0.00, 1.9 ) -- (0.00, 2.1);
\draw[] (2.7, 1.9 ) -- (2.7, 2.1);

\node[below, text=black] at (-1.75, 1.9) {$\sigma_1$};
\node[below, text=black] at (0, 1.9)     {$\sigma_2$};
\node[below, text=black] at (2.7, 1.9)   {$\sigma_3$};

\addplot[only marks] table [
  x expr = {\thisrowno{0}},
  y expr = {2.1}
] {
  -3
  -2
  -1.5
  -1

  2
  2.2
  2.4
  2.6
  2.8
  3
  3.2
  3.4
};

\draw[] (-1.875, 0.9 ) -- (-1.875, 1.1);
\draw[] (2.3, 0.9 ) -- (2.3, 1.1);
\draw[] (3.1, 0.9 ) -- (3.1, 1.1);

\node[below, text=black] at (-1.875, 0.9) {$\sigma_1$};
\node[below, text=black] at (2.3, 0.9)     {$\sigma_2$};
\node[below, text=black] at (3.1, 0.9)   {$\sigma_3$};

\addplot[only marks] table [
  x expr = {\thisrowno{0}},
  y expr = {1.1}
] {
  -2
  -1.75
  -1.5
  -1.25

  1.4
  1.6
  1.8
  2.0
  2.2
  2.4
  2.6
  2.8
};

\draw[] (-1.625, -0.1 ) -- (-1.625, 0.1);
\draw[] (1.7, -0.1 ) -- (1.7, 0.1);
\draw[] (2.5, -0.1 ) -- (2.5, 0.1);

\node[below, text=black] at (-1.625, -0.1) {$\sigma_1$};
\node[below, text=black] at (1.7, -0.1)     {$\sigma_2$};
\node[below, text=black] at (2.5, -0.1)   {$\sigma_3$};

\addplot[only marks] table [
  x expr = {\thisrowno{0}},
  y expr = {0.1}
] {
  -2
  -1.75
  -1.5
  -1.25

  1.6
  1.8
  2
  2.2
  2.4
  2.6
  2.8
  3
};

\draw[] (-1.625, -1.1 ) -- (-1.625, -0.9);
\draw[] (1.9, -1.1 ) -- (1.9, -0.9);
\draw[] (2.7, -1.1 ) -- (2.7, -0.9);

\node[below, text=black] at (-1.625, -1.1) {$\sigma_1$};
\node[below, text=black] at (1.9, -1.1)     {$\sigma_2$};
\node[below, text=black] at (2.7, -1.1)   {$\sigma_3$};

\addplot[only marks] table [
  x expr = {\thisrowno{0}},
  y expr = {-0.9}
] {
  -2
  -1.75
  -1.5
  -1.25

  1.6
  1.8
  2
  2.2
  2.4
  2.6
  2.8
  3
};

\end{axis}
\end{tikzpicture}
\end{minipage}
\caption{A graphical representation of the shift migration process throughout the SCF
procedure. The SCF iterations progress from top to bottom with the filled circles
representing the validated Ritz values at that iteration. $\sigma_{1}$, $\sigma_2$, $\sigma_3$ represent
the spectral shifts used to obtain the Ritz values at each iteration. The shifts
for the first SCF iteration are representative of a typical DOS-shift selection
scheme where shifts are chosen both in regions without eigenvalues as well
as regions with a disproportionately large number of eigenvalues due to inaccuracy
in the DOS approximation. At each subsequent SCF iteration, new shifts are chosen
via k-means clustering of the Ritz pairs obtained from the previous iteration.}
\label{fig:shift_migration}
\end{figure}
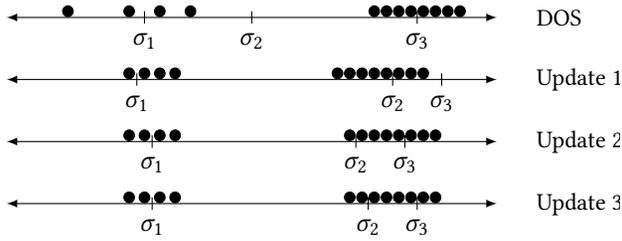

Due to the limited number of Ritz values that can be extracted from the Lanczos
method in the spectral interior, it is possible that the initial selection of
spectral shifts produced by a Lanczos DOS estimation procedure is far from
optimal.  In particular, shifts may be placed in spectral regions devoid of
eigenvalues. Another possible scenario is that an insufficient number of shifts
are placed in regions that contain a disproportionately large number of
eigenvalues.  An illustration of this issue is given on the axis labeled ``DOS"
in \cref{fig:shift_migration}.  However, shift misplacement can be
incrementally corrected in subsequent SCF iterations by using a clustering
algorithm to partition previously computed eigenvalue approximations and refine
the shift selection.


For each SCF iteration, we obtain a set of eigenpair approximations for $(\Am{i},B)$.
Thus, for $i>0$, we have available to us a set of approximate eigenvalues for $(\Am{i-1},B)$.
In the early SCF iterations, when the change in eigensystem between two subsequent iterations
is relatively large, it is possible that the shifts selected for $(\Am{i-1},B)$ would not
be appropriate for the slicing of the spectrum for $(\Am{i},B)$. We discuss a strategy to determine
this (dis)similarity and strategies for subsequent shift selection in \cref{sec:misev}.
However, if the eigenvalues of $(\Am{i},B)$ are sufficiently close
to those of $(\Am{i-1},B)$, then we may use the approximated eigenvalues of $(\Am{i-1},B)$
as a reference to determine the spectral shift placement for the spectrum slicing of
$(\Am{i},B)$.
Due to the localized nature of the shift-invert spectral transformation, 
rapid convergence of the shift-invert subspace iteration is achieved when
spectral shifts are placed centrally in clusters of eigenvalues.
Thus, we may determine more optimal shift placement by identifying spectral clusters from
the computed eigenvalues $(\Am{i-1},B)$ and placing shifts in the centroids
of these clusters for determination of the eigenpairs of $(\Am{i},B)$.

To identify spectral clusters, we employ the k-means clustering algorithm \cite{lloyd1982least}.
At the $i$-th SCF iteration for $i>0$, we use k-means to
identify $n_s$ clusters from the validated eigenpairs obtained from the $(i-1)$-st iteration.
The centroids of the clusters may then be used in the generation of the
$i$-th set of spectral shifts.
As the SCF procedure converges, the centroids  of the clusters will also
converge to a particular set of spectral shifts. An illustration of this
convergence behavior is given in \cref{fig:shift_migration}.


Although the k-means clustering problem is generally NP-hard,
we do not necessarily need to obtain a globally optimal solution
to the clustering problem in order to identify appropriate
spectral shifts.
Our objectives are to identify eigenvalue clusters and to partition
 nearly uniformly distributed eigenvalues into slices of roughly equal size.
In general, determination of $n_s$ clusters is a drastic over clustering of the
Ritz values. However, k-means clustering usually results in equal sized
clusters, even in the case of over clustering.
Due to the fact that $n_s$ is relatively small, obtaining clusters from
this data using k-means may be achieved with negligible cost.

The k-means algorithm is an iterative procedure initialized with
a set of guesses to cluster centroids.  The choice of these initial guesses can
have a significant effect on the convergence of the algorithm and the quality
of cluster centroids it produces.  In the SISLICE method, these guesses are
usually taken to be the spectral shifts used in the previous SCF iteration.
However, if the previous spectral shifts are generated from the DOS shift
selection strategy, it is possible that the k-means algorithm can converge to a
sub-optimal solution if the centroids are initialized with these shifts.
To address this issue, we employ the k-means++ \cite{arthur2007kmeans} cluster
initialization strategy to improve initial guesses of the centroids prior to
the k-means clustering process. Rather than select all guess centroids at
random with uniform probability for all eigenvalues, k-means selects the first
centroid with a uniform probability and then selects all subsequent centroids
from the remaining eigenvalues according to a probability distribution which is
quadratic in the distance of each eigenvalue to its nearest existing guess
centroid. This process is repeated until $n_s$ centroids have been selected.
The result is a set of guess centroids with are well separated and are selected
such that there is a provable upper bound to the k-means objective function.
For more details, we refer the reader to \cite{arthur2007kmeans}.

Assuming that the eigenvalue distribution does not change drastically
throughout the SCF procedure, k-means also allows the SISLICE method to track 
eigenvalue changes and migrate shifts
between SCF iterations as shown in the axes labeled ``Update" in
\cref{fig:shift_migration}. Until convergence is reached, shift migration is
performed by k-means using the validated eigenpairs of the previous SCF iteration
as the initial guess centroids.
As k-means converges to the local optimum nearest to the initial guess,
this choice of guess allows for the k-means shift migration strategy to
converge to a single set of shifts as the SCF converges.
It is
important that the clustering is performed on validated eigenpairs to avoid
oversampling of spectral regions for which the Ritz pairs of adjacent spectral
probes overlap or contain spurious Ritz values. We demonstrate the efficacy of
this migration scheme in \cref{sec:numerical}.
%

Once a clustering of the validated Ritz values has been obtained, we may use
the centroids of those clusters to generate the set of spectral shifts for
the next SCF iteration.
Instead of creating new spectral probes, we would like
to reuse the Ritz vectors produced by the existing spectral probes as
the initial guesses to the desired eigenvectors in the subsequent
shift-invert subspace iteration to improve convergence. Therefore, in the SISLICE method,
we update the shifts of the existing spectral probes based on the clustering
information rather than starting completely from scratch.
To update the shift for each spectral probe,
we form a mapping between the eigenvalue clusters and spectral
probes such that each cluster is mapped to the probe with which
it has maximal overlap.

%

In the case that this map is bijective, the spectral shift associated with the
spectral probe is taken to be the centroid of its associated cluster.  However,
it is possible, especially in the early SCF iterations when the
eigenspectrum undergoes considerable change, that this map is not bijective.
As a result, there is some ambiguity as to how to best update the probes
which have no preimage under this mapping.

\begin{figure}[t]
  \centering
  \begin{minipage}{0.49\textwidth}
  \centering
  \begin{tikzpicture} 
\begin{axis}[ 
    axis x line=none, 
    axis y line=none, 
    ymin=-4, ymax=4,
    xmin=-5, xmax=5,
    axis equal image,
    xticklabels={},
    every tick/.style={
        black,
        thick,
    },
    width=\textwidth,
    height=\textwidth
]

\draw[>=latex, <->] (-5,2) -- (5,2);
\draw[dashed]       (0, 1.5) -- (0, 2.5);
\draw[]  (-2.5,1.8) -- (-2.5,2.2);
\draw[]  (2.5,1.8) -- (2.5,2.2);
\node[below,text=black] at (-3,1.8) {$\sigma_1$};
\node[below,text=black] at (3,1.8) {$\sigma_2$};
\addplot[only marks,color=blue] table [
  x expr = {\thisrowno{0}},
  y expr = {2.125}
] {
  -0.50
  0.25
  0.75
  1.25
};
\addplot[only marks,color=red] table [
  x expr = {\thisrowno{0}},
  y expr = {2.125}
] {
2.75
3.25
3.75
4.25
};

\draw[>=latex, <->] (-5,-2) -- (5,-2);
\draw[dashed]       (1.96875, -2.5) -- (1.96875, -1.5);
\draw[]  (0.4375,-2.2) -- (0.4375,-1.8);
\draw[]  (3.5,-2.2) -- (3.5,-1.8);
\node[below,text=black] at (0.625,-2.2) {$\sigma_1$};
\node[below,text=black] at (3.5,-2.2) {$\sigma_2$};
\addplot[only marks,color=blue] table [
  x expr = {\thisrowno{0}},
  y expr = {-1.875}
] {
  -0.50
  0.25
  0.75
  1.25
};
\addplot[only marks,color=red] table [
  x expr = {\thisrowno{0}},
  y expr = {-1.875}
] {
2.75
3.25
3.75
4.25
};
\end{axis} 
\end{tikzpicture}
  \end{minipage}~
  \begin{minipage}{0.49\textwidth}
  \caption{A schematic illustration of how multiple clusters may be mapped to
    the same probe, and how an old probe may be deleted and a new probe can be
    inserted.  The blue and red dots are two clusters of approximate eigenvalues
    that are both mapped to the same spectral probe centered at $\sigma_2$.  After
    the mapping between clusters and the previous spectral probes is established,
    the probe centered at $\sigma_1$ is deleted because no cluster is mapped to it.
    A new problem centered at a new shift $\sigma_1$ is inserted, and $\sigma_2$ is
    also moved to the right.}
  \label{fig:clustermap}
  \end{minipage}
\end{figure}
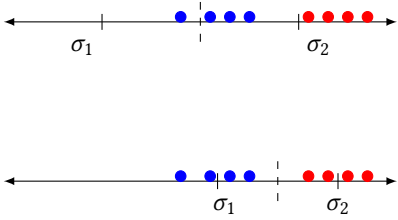

The case of non-bijective maps between clusters and spectral probes is
typically a symptom of poor shift selection in the previous SCF iteration
resulting in some spectral probes picking up only a few validated Ritz values while others capturing a disproportionate number of validated Ritz values. As the validated Ritz values are
separated by the k-means algorithm into different clusters, validated Ritz values retrieved from
one spectral probe may be separated into several clusters resulting in
several clusters being mapped to the same spectral probe (see the
schematic illustration in \cref{fig:clustermap}.)
In the meantime, the few validated Ritz values obtained by a poorly
placed spectral probe SP($\sigma_j$) may be placed into a cluster
that gets mapped to a different spectral probe SP($\sigma_j$), leaving
SP($\sigma_j$) without any cluster to map to.

When several clusters are mapped to the same spectral probe, we merge
these clusters into a single cluster.  The resulting centroid of the merged
Ritz values is taken to be the spectral shift of the associated spectral probe. In
the case that the mapped cluster contains too many eigenvalues, a probe may
be inserted to ensure proper load balance.

If a spectral probe is not mapped to any cluster, it is
simply deleted.  However, since the number of spectral shifts (and thus the number of spectral probes) is fixed throughout the SCF procedure,
when a probe is removed, another probe must be added to maintain this
fixed number of shifts. We choose to add a probe to a cluster that
contains the largest number of Ritz values.
To add such a spectral probe to such a cluster, we need to break up the cluster
first into two clusters. In practice,
this may be achieved by performing a 2-means clustering of the Ritz pairs in the
largest cluster. One of the clusters is mapped to the original spectral
probe mapped to the cluster before it was broken up. The shift associated
with that cluster is replaced with the centroid of the new smaller cluster.
The other cluster is mapped to the newly added spectral probe. In addition
to setting the shift of the probe to the centroid of the new cluster, we
also need to copy (or send in a distributed memory implementation) the Ritz vectors associated with the Ritz values in
this cluster to the added spectral probe.
This process of breaking up a large cluster and adding a new spectral probe
is repeated until the desired number of spectral
shifts and probes is obtained. The specifics regarding probe insertion in a
distributed infrastructure are discussed \cref{sec:parallel}.

The SIESTA-SIPs method of \cite{keceli18_1806} adopts a similar k-means
strategy for tracking the eigenvalue migration throughout the SCF procedure to
avoid costly methods to perform shift selection every SCF iteration.  The
SIESTA-SIPs method is not as sensitive to changes in the definitions of
spectral slice intervals between SCF iterations as long as they are distributed
such to allow for similar convergence rates between slices due to their use of
SI-Lanczos.  This also inhibits their ability to leverage the knowledge of the
convergent nature of the eigenvectors to accelerate their eigenpair convergence
as the SCF iterations progress.  This leads to two main differences in the
utilization of k-means between SISLICE and SIESTA-SIPs:
\begin{enumerate}
  \item Whereas SISLICE aims to minimize the change in centroid location between
  SCF iterations by initializing the centroids as the shifts from the previous
  SCF iteration, SIESTA-SIPs obtains a new and possibly different set of centroids
  every iteration based on the convergence of k-means using a uniform initial guess.

  \item The fact that they do no utilize previous eigenvector information means
  that they need not map the newly obtained cluster centroids to existing sets of
  spectral probes as this information is not generally passed between SCF
  iterations in their prescription. 
\end{enumerate}

\subsection{Missing Eigenvalues}
\label{sec:misev}
Because $A^{(i)}$ can change significantly from $A^{(i-1)}$ in early SCF
iterations, a shift selection scheme based on the clustering of approximate
eigenvalues
 of $A^{(i-1)}$ may not be optimal for computing eigenpairs of $A^{(i)}$.
In particular, it is possible that
spectral probes constructed from the sub-optimal selection of target shifts
miss some eigenvalues.  The spectral slices in which these
missing eigenvalues reside can be identified in the validation
process described in \cref{sec:valid}.

When missing eigenvalues are detected, we perform a new DOS estimation on
$(A^{(i)},B)$ with an appropriate resolution to place new shifts in
spectral slices that contain missing eigenvalues.  New spectral probes
are created to recompute approximate eigenvalues within these newly created
spectral slices.  This is a costly step because the next SCF cycle cannot
start until all missing eigenvalues are accounted for. An example of this
state of affairs is demonstrated in \cref{sec:misev_results}.



To reduce the likelihood of missing eigenvalues resulting from shift
misplacement, we can monitor the convergence of SCF for drastic changes
in the spectrum by comparing the partial traces of the system matrices
within the subspace spanned by the previously validated eigenvectors.
For example, given the metric
\begin{equation}
\eta(V,A) = \mathrm{Tr}\left( V^{T} A V \right),
\end{equation}
if the difference between $\eta(X^{(i-1)},\Am{i-1})$ and $\eta(X^{(i-1)},\Am{i})$ is larger than
some specified threshold, then the spectra of the matrices may be deemed to be sufficiently
dissimilar. If this is found to be the case, then it would be beneficial to use the DOS
shift selection strategy discussed in \cref{sec:dos}. Although this strategy
does not completely eliminate the possibility of missing eigenvalues
(because individual eigenvalues can move around without affecting
the trace of the $A^{(i)}$), it may help reduce that possibility
and the cost associated with generating new probes to seek the missing
eigenvalues.

\section{Parallel Implementation}
\label{sec:parallel}

The algorithmic subtasks described in the previous sections have been constructed
in such a way as to allow for maximal concurrency in the slicing of the
spectral region of interest: each of the spectral probes may be constructed
independently of any other spectral probe.  As the construction of the spectral
probes through shift-invert subspace iterations constitutes the bulk of the
work in the SISLICE method, this task independence should lead to
scalable performance.  The slice validation
scheme outlined in \cref{sec:valid} would require some level of synchronization
between independent computing units. Further, the shift insertion and deletion
schemes outlined in \cref{sec:kmeans} would require some data to be copied
from some processors/nodes to others. However, these communication and
synchronization overheads are generally small as we will see in the next section.
In this section, we outline the salient aspects of the parallel implementation of
the SISLICE method.

We note that the parallelization discussed here focuses exclusively
on the parallel execution of spectral probes. An additional
level of finer grain parallelism exists within each spectral probe.
If matrices $A$ and $B$ can be replicated and stored on each single many-core
compute node, a hybrid-parallelism scheme utilizing both shared-memory and
message passing parallelism may be achieved through exploitation of optimized
implementations of threaded BLAS and LAPACK (such as those found in Intel(R)
MKL, IBM(R) ESSL, OpenBLAS, BLIS, ATLAS, cuBLAS, etc) within a particular MPI rank.
While we do not treat this level of parallelism explicitly in this section, its
leverage is trivial on modern computing architectures and is implied for the numerical experiments
in \cref{sec:numerical}. If $A$ and $B$ are too large to be stored on
a single compute node, then the factorizations and linear system solves required
for each probe may be performed using ScaLAPACK in the case of dense matrices,
or a distributed sparse solver such as symPACK, MUMPS,
or PARDISO in the case of sparse matrices. We should note that
the parallel scalability of $LDL^T$ factorization and back substitutions
for solving triangular systems with multiple right hand sides is
generally much better than what can be achieved in a dense eigensolver.
These distributed calculations for each spectral probe may take place on a
subset of the total number of MPI ranks, allowing leverage of massive
parallelism on large computing clusters.  We do not treat this particular
parallelism scheme in this work, but it has been discussed at length in other
related work \cite{zhang2007sips,keceli16_448,keceli18_1806}.

\subsection{Spectral Probe Distribution and Synchronization}
\label{sec:parallel_sync}

The SISLICE method is designed for taking advantage of computer systems that
have a large amount of computational resources in terms of compute nodes
and cores within each node.  In an ideal scenario, the number of spectral
probes should match the number of nodes (or groups of nodes) so that
all probes can be executed simultaneously.  
The optimal number of nodes corresponding to the natural clustering the eigenvalues
could be determined by the DOS based shift selection in \cref{sec:dos} and queried
prior to running the SISLICE solver. However, we note that this is not a requirement
of the implementation discussed in this work.

When the number of computational nodes is less than the number of
spectral probes, a round-robin distribution of probes to nodes can be
used, i.e., SP($\sigma_j$) may be mapped to the ($j$ mod $n_r$)-th MPI rank,
where $n_r$ is the number of MPI ranks.  In this case, the computation
is not load balanced if $n_r$ does not divide $n_s$.

Once the spectral probes have been constructed, each SP($\sigma_j$) contains a
set of Ritz pairs which approximate the eigenpairs in the neighborhood of
$\sigma_j$. However, the slice validation scheme described in \cref{sec:valid}
requires knowledge of Ritz pair information from adjacent probes.
If each of the adjacent
spectral probes has been constructed on a different MPI rank (or group of MPI
ranks), the validation scheme requires some level of communication /
synchronization of Ritz pair information between the MPI ranks.
However, the validation scheme only requires knowledge of the Ritz values and
associated residual norms to validate the spectral slices. The Ritz vectors are
not explicitly required.

If the Ritz values and residual norms were only to be used in the slice
validation scheme, their synchronization could be further limited to only the
neighboring ranks of the owner of a particular spectral probe.  However,
because the entire set of validated Ritz values is used in updating spectral
shifts through k-means clustering (as described in \cref{sec:kmeans}), it is
useful to synchronize this information across all of the MPI ranks. While a
distributed implementation of k-means clustering is possible, the fact that
each spectral probe only accounts for a relatively small number of validated
Ritz pairs would require excessive communication to perform the clustering.
Because the storage requirement of the Ritz values and residual norms is
negligible relative to the Ritz vectors, this synchronization scheme poses no
storage overhead.  In the SISLICE method, the synchronization of Ritz values
and residual norms for each spectral probe may be implemented as
\texttt{MPI\_Allgather}, and poses only minimal communication overhead relative
to the computational cost of the shift-invert subspace iterations. As the Ritz
values and residual norms are replicated across each MPI rank, the tasks of
slice validation and shift updates may also be replicated to avoid
communication. We note for clarity that in the case of random initialization of
the clustering problem through e.g. k-means++, the Ritz value clustering may
still be replicated through the use of pseudo random number generation using
the same seed value. The scalability of this distribution and synchronization
scheme is demonstrated in \cref{sec:scaling}.


\subsection{Spectral Probe Insertion and Removal}
\label{sec:probe_insert}

As was discussed in \cref{sec:kmeans}, occasionally
shift selection and migration schemes employed by the SISLICE method yield
sub-optimal shift placement leading to redundant spectral probes and probes
which are responsible for a disproportionate number of validated eigenpairs.
This is typically the case in the early SCF iterations due to
the crude DOS approximation by the Lanczos procedure described in \cref{sec:dos}.
The presence of redundancies in the spectral probes leads to a load imbalance
which should be avoided to ensure scalability on large, distributed  computing
systems.  For the purposes of this section, the term ``load balance" should be
thought of as balance of \emph{useful} work. Technically speaking, even in the
case of redundancies in the spectral probes, the computational work performed
for each spectral probe will always be roughly the same given that the number of
subspace iterations performed and subspace dimensions are uniform across all probes. Thus
this work is always ``balanced". However, we want to ensure that each rank is
performing a roughly equal amount of useful work (in the sense of yielding a
roughly equal number of validated Ritz pairs) rather than wasting valuable
computational resources in spectral regions where it is not needed.

As the number of spectral probes is fixed in the SISLICE method, removal
of a spectral probe necessitates the insertion of a spectral probe to
balance the work in another spectral region. This probe removal necessarily
leads to a load imbalance if the work was balanced in the previous SCF iteration.
As discussed in \cref{sec:kmeans},
probes are inserted so as to break up large clusters of validated Ritz values
such that they effectively span multiple probes after the subsequent shift-invert subspace iterations.
In a distributed computing environment, care must be taken
to ensure probe insertion is performed in such a way as to balance the work
between independent computing ranks while avoiding a large communication overhead.
As the k-means clustering is replicated on each rank, the decision of where to insert new spectral probes may also
be effectively replicated with only minimal communication. For each probe to be
inserted, the determination of the two new spectral shifts is replicated on
each rank. The new probe which is to be inserted is assigned to the rank with
the least amount of work (thus ensuring load balance). Once this decision has been made,
the probe whose shift has been moved through this procedure communicates its
Ritz vector data to the newly inserted probe to allow its reuse in the subsequent
SCF iteration. The cost of this point-to-point communication is relatively small
in practice and may be overlapped with
the determination and communication of other probe insertions.

In the case when  shifts must be inserted due to missing eigenvalues within a particular
spectral slice as described in \cref{sec:misev}, we may leverage the fact that the computation
is done in parallel to our advantage. It may be the case that the DOS shift insertion strategy
yielded several shifts in the spectral region that contains the missing eigenvalues. Rather than
have processors or processor groups sit idle while the missing eigenpairs are obtained
sequentially, the newly inserted probes may be distributed in the same manner as the initial
probe distribution. Due to the fact that missing eigenvalues are typically a symptom of poorly
placed shifts, not of too few shifts, inserting probes will not yield
more useful probes than $n_s$, i.e. if a probe had to be inserted to resolve missing eigenpairs,
it is typically the case that some probes did not produce validated eigenpairs. However, even
if each of the probes from the first round of subspace iterations produced validated eigenpairs,
the mapping scheme between eigenvalue clusters and spectral probes will preclude the possibility
of yielding more than $n_s$ probes for the subsequent SCF iteration. This is due to the fact that
the SISLICE method obtains $n_s$ clusters regardless of the number of probes used to produce
the validated eigenpairs in the previous SCF iteration. Thus, even in the case of probe insertion
in the previous SCF iteration, the SISLICE method ensures load balance is maintained in subsequent SCF iterations.

\section{Numerical Experiments}
\label{sec:numerical}

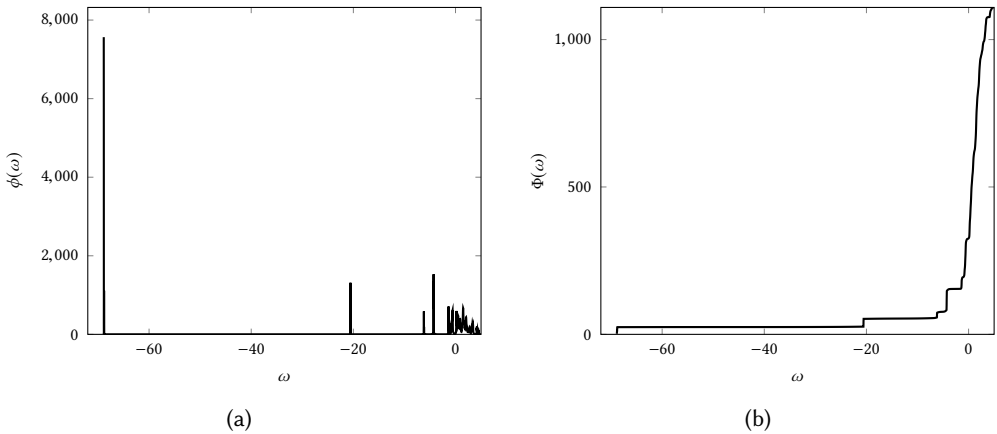
\begin{figure}[t]

  \centering
  \begin{subfigure}[b]{0.49\textwidth}
    \centering
    \begin{tikzpicture}
\pgfmathsetlengthmacro\MajorTickLength{
  \pgfkeysvalueof{/pgfplots/major tick length} * 0.5
}
\begin{axis}[
  width=\textwidth,
  every axis plot/.append style={ thick},
  legend style={font=\tiny},
  legend cell align={left},
  every tick label/.append style={font=\tiny},
  major tick length=\MajorTickLength,
  xmin=-72, xmax=5,
  ymin = 0,
  xlabel={$\omega$},
  ylabel={$\phi(\omega)$},
  x label style={ font=\tiny },
  y label style={ font=\tiny },
]
\addplot[mark={}] table [x=x, y=y, col sep=comma] {figs/SiOSi5_fixed_dos.dat};
\end{axis}
\end{tikzpicture}
    \caption{}
    \label{fig:silane_dos}
  \end{subfigure}
  \begin{subfigure}[b]{0.49\textwidth}
    \centering
    \begin{tikzpicture}
\pgfmathsetlengthmacro\MajorTickLength{
  \pgfkeysvalueof{/pgfplots/major tick length} * 0.5
}
\begin{axis}[
  width=\textwidth,
  every axis plot/.append style={ thick},
  legend style={font=\tiny},
  legend cell align={left},
  every tick label/.append style={font=\tiny},
  major tick length=\MajorTickLength,
  xmin=-72, xmax=5,
  ymin = 0, ymax = 1109,
  xlabel={$\omega$},
  ylabel={$\Phi(\omega)$},
  x label style={ font=\tiny },
  y label style={ font=\tiny },
  ylabel shift = -4 pt,
]
\addplot[mark={}] table [x=x, y=y, col sep=comma] {figs/SiOSi5_fixed_cdos.dat};
\end{axis}
\end{tikzpicture}
    \caption{}
    \label{fig:silane_cdos}
  \end{subfigure}

  \caption{Lanczos DOS (a) and CDOS (b) for the entire spectrum of Silane ($\probdim=1109$).
    Silane exhibits a number of isolated eigenvalue clusters (spikes in the DOS) lower
    in the spectrum and a more uniform distribution at larger eigenvalues. The DOS and
    CDOS calculations were performed using 100 Lanczos iterations with the converged matrix pencil.
  }
  \label{fig:silane}

\end{figure}

\begin{figure}[t]

  \centering
  \begin{subfigure}[b]{0.49\textwidth}
    \centering
    \begin{tikzpicture}
\pgfmathsetlengthmacro\MajorTickLength{
  \pgfkeysvalueof{/pgfplots/major tick length} * 0.5
}
\begin{axis}[
  width=\textwidth,
  every axis plot/.append style={ thick},
  legend style={font=\tiny},
  legend cell align={left},
  every tick label/.append style={font=\tiny},
  major tick length=\MajorTickLength,
  xmin=-1.73, xmax=-0.64,
  ymin = 0,
  xlabel={$\omega$},
  ylabel={$\phi(\omega)$},
  x label style={ font=\tiny },
  y label style={ font=\tiny },
]
\addplot[mark={}] table [x=x, y=y, col sep=comma] {figs/Graphene720_fixed_dos.dat};
\end{axis}
\end{tikzpicture}
    \caption{}
    \label{fig:graphene_dos}
  \end{subfigure}
  \begin{subfigure}[b]{0.49\textwidth}
    \centering
    \begin{tikzpicture}
\pgfmathsetlengthmacro\MajorTickLength{
  \pgfkeysvalueof{/pgfplots/major tick length} * 0.5
}
\begin{axis}[
  width=\textwidth,
  every axis plot/.append style={ thick},
  legend style={font=\tiny},
  legend cell align={left},
  every tick label/.append style={font=\tiny},
  major tick length=\MajorTickLength,
  xmin=-1.73, xmax=-0.64,
  ymin = 0, ymax = 1109,
  xlabel={$\omega$},
  ylabel={$\Phi(\omega)$},
  x label style={ font=\tiny },
  y label style={ font=\tiny },
  ylabel shift = -4 pt,
]
\addplot[mark={}] table [x=x, y=y, col sep=comma] {figs/Graphene720_fixed_cdos.dat};
\end{axis}
\end{tikzpicture}
    \caption{}
    \label{fig:graphene_cdos}
  \end{subfigure}

  \caption{Lanczos DOS (a) and CDOS (b) for the lowest 1000 eigenvalues of Graphene ($\probdim=9360$).
    Graphene exhibits a nearly uniform eigenvalue distribution throughout its entire spectrum.
    The DOS and CDOS calculations were performed using 100 Lanczos iterations with the converged matrix pencil.
  }
  \label{fig:graphene}

\end{figure}
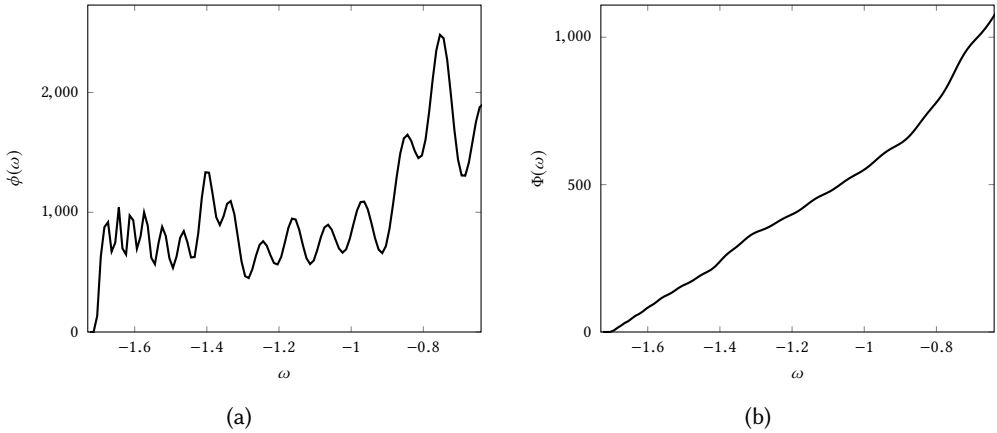

In this section, we report a set of numerical experiments which demonstrate the effectiveness of the proposed shift selection technique
for computing all or a subset of eigenvalues of a
matrix pencil or a sequence of matrix pencils.
We examine two limiting cases of eigenvalue distribution shown in
\cref{fig:silane,fig:graphene}.

The Silane test case (\cref{fig:silane}, $\probdim=1109$) is an all-electron density functional
theory calculation using a Gaussian basis set.  Its spectrum exhibits a number
of isolated eigenvalue clusters at lower eigenvalues and a more uniform
distribution at larger eigenvalues. The isolated eigenvalue clusters at low
eigenvalues are a common feature in all-electron density functional
calculations. All matrices related to the Silane test case in this work were
obtained using the NWChemEx software package \cite{nwchemex}.

The Graphene test case (\cref{fig:graphene}, $\probdim=9360$) is a density functional theory
calculation using pseudo-potentials for the core electrons. As such, its
spectrum does not contain isolated clusters. The more uniform nature of the
spectrum is a
common feature in nearly all pseudo-potential based density functional calculations.
All matrices related to the Graphene test case in this work were obtained using
the SIESTA software package \cite{soler2002siesta}.

\subsection{Shift Selection for a Fixed Matrix Pencil}
\label{sec:fixed}

To demonstrate how the shift selection strategy enables rapid convergence
of the SISLICE method as the SCF procedure approaches convergence,
i.e. when the matrix pencils change very
little between SCF iterations, we examine the
convergence of eigenpairs for fixed matrix pencils in this section.
This experiment allows us to gauge the effectiveness of the
shift selection strategy when more accurate estimation of
the desired eigenvalues becomes available in successive
SCF iterations.

To simplify our exposition, we examine a set of representative spectral windows for the
aforementioned test cases using the converged matrix pencil $(A,B)$ as the
representative eigenvalue problem.  Even though the matrix pencil
does not change, an artificial SCF procedure is carried out
and a new set of shifts may be chosen after a fixed number of
subspace iterations have been performed.
This procedure can be viewed as a generalized (block) Rayleigh quotient
iteration.

All calculations in this section were
performed using a probe basis dimension of $k = 100$ and 4 shift-invert subspace iterations
per SCF iteration.  The SCF iteration is
considered converged if the maximum of the residual norms associated with
all validated approximate eigenpairs is below the threshold of $10^{-13}$.
In both of the presented cases, we
observe the expected monotonic convergence of the eigenvalues within individual
spectral slices once shift migration has been performed.

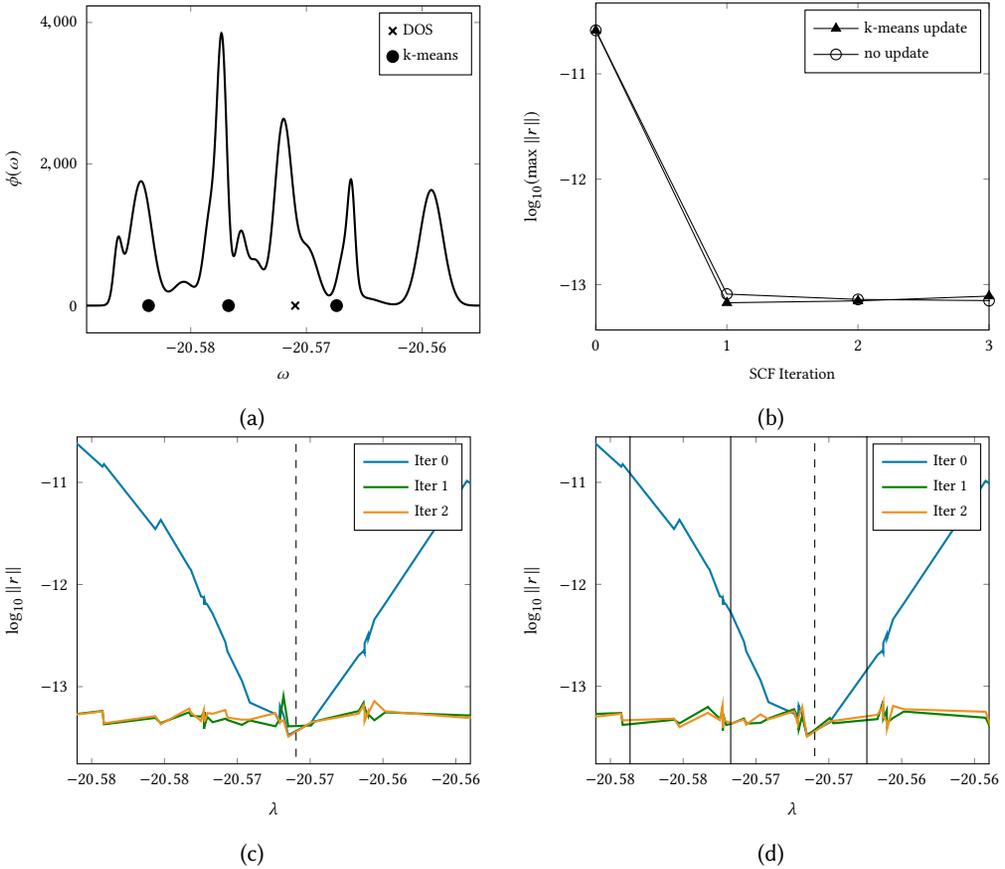
\begin{figure}[t]

  \begin{subfigure}[b]{0.49\textwidth}
    \begin{tikzpicture}
\pgfmathsetlengthmacro\MajorTickLength{
  \pgfkeysvalueof{/pgfplots/major tick length} * 0.5
}
\begin{axis}[
  width=\textwidth,
  every axis plot/.append style={ thick},
  legend style={font=\tiny},
  legend cell align={left},
  every tick label/.append style={font=\tiny},
  major tick length=\MajorTickLength,
  xmin=-20.589, xmax=-20.555,
  ylabel={$\phi(\omega)$},
  xlabel={$\omega$},
  x label style={ font=\tiny },
  y label style={ font=\tiny }
]
\addplot[only marks, mark=x] table {
-20.570965906 0
};
\addplot[only marks, mark=*] table {
-20.583646653 0
-20.576737164 0
-20.567389376 0
};
\legend{
DOS,
k-means
};
\addplot[mark={}] table [x=x, y=y, col sep=comma] {figs/SiOSi5_fixed_core_cluster_dos.dat};
\end{axis}
\end{tikzpicture}
    \caption{}
    \label{fig:silane_core_cluster_fixed_dos}
  \end{subfigure}~
  \begin{subfigure}[b]{0.49\textwidth}
    \begin{tikzpicture}
\pgfmathsetlengthmacro\MajorTickLength{
  \pgfkeysvalueof{/pgfplots/major tick length} * 0.5
}
\begin{axis}[
  width=\textwidth,
  legend style={font=\tiny},
  legend cell align={left},
  every tick label/.append style={font=\tiny},
  major tick length=\MajorTickLength,
  xmin=0, xmax=3,
  ylabel={$\log_{10}(\mathrm{max}\,\|r\|)$},
  xlabel={SCF Iteration},
  xtick={0,1,2,3},
  x label style={ font=\tiny },
  y label style={ font=\tiny }
]
\addplot[mark=triangle*,fill opacity=1] table [x=x, y=y1, col sep=comma] {figs/SiOSi5_fixed_waterfall.dat};
\addplot[mark=o,fill opacity=0] table [x=x, y=y2, col sep=comma] {figs/SiOSi5_fixed_waterfall.dat};
\legend{
k-means update, 
no update,
};
\end{axis}
\end{tikzpicture}
    \caption{}
    \label{fig:silane_core_cluster_fixed_waterfall}
  \end{subfigure}~\\

  \begin{subfigure}[b]{0.49\textwidth}
    \begin{tikzpicture}
\pgfmathsetlengthmacro\MajorTickLength{
  \pgfkeysvalueof{/pgfplots/major tick length} * 0.5
}
\begin{axis}[
  width=\textwidth,
  legend style={font=\tiny},
  legend cell align={left},
  every axis plot/.append style={ thick},
  every tick label/.append style={font=\tiny},
  major tick length=\MajorTickLength,
  xmin=-20.586, xmax=-20.559,
  ylabel={$\log_{10}\|r\|$},
  xlabel={$\lambda$},
  x label style={ font=\tiny },
  y label style={ font=\tiny }
]
\addplot[mark={},color=cerulean] table [x=x, y=y1, col sep=comma] {figs/SiOSi5_fixed_core_cluster_conv_nokmeans.dat};
\addplot[mark={},color=applegreen] table [x=x, y=y2, col sep=comma] {figs/SiOSi5_fixed_core_cluster_conv_nokmeans.dat};
\addplot[mark={},color=carrot] table [x=x, y=y3, col sep=comma] {figs/SiOSi5_fixed_core_cluster_conv_nokmeans.dat};
\draw[dashed] (-20.570965906, -14) -- (-20.570965906, -10);
\legend{
Iter 0,
Iter 1,
Iter 2
};
\end{axis}
\end{tikzpicture}
    \caption{}
    \label{fig:silane_core_cluster_fixed_conv_nokmeans}
  \end{subfigure}~
  \begin{subfigure}[b]{0.49\textwidth}
    \begin{tikzpicture}
\pgfmathsetlengthmacro\MajorTickLength{
  \pgfkeysvalueof{/pgfplots/major tick length} * 0.5
}
\begin{axis}[
  width=\textwidth,
  legend style={font=\tiny},
  legend cell align={left},
  every axis plot/.append style={ thick},
  every tick label/.append style={font=\tiny},
  major tick length=\MajorTickLength,
  xmin=-20.586, xmax=-20.559,
  ylabel={$\log_{10}\|r\|$},
  xlabel={$\lambda$},
  x label style={ font=\tiny },
  y label style={ font=\tiny }
]
\addplot[mark={},color=cerulean] table [x=x, y=y1, col sep=comma] {figs/SiOSi5_fixed_core_cluster_conv.dat};
\addplot[mark={},color=applegreen] table [x=x, y=y2, col sep=comma] {figs/SiOSi5_fixed_core_cluster_conv.dat};
\addplot[mark={},color=carrot] table [x=x, y=y3, col sep=comma] {figs/SiOSi5_fixed_core_cluster_conv.dat};
\draw[dashed] (-20.570965906, -14) -- (-20.570965906, -10);
\draw[] (-20.583646653, -14) -- (-20.583646653, -10);
\draw[] (-20.576737164, -14) -- (-20.576737164, -10);
\draw[] (-20.567389376, -14) -- (-20.567389376, -10);
\legend{
Iter 0,
Iter 1,
Iter 2
};
\end{axis}
\end{tikzpicture}
    \caption{}
    \label{fig:silane_core_cluster_fixed_conv}
  \end{subfigure}

  \caption{Isolated cluster of Silane eigenvalues $[-20.59, -20.55]$ (37 eigenvalues). (a) Initial Lanczos DOS along with
    DOS shift placement and k-means update. (b) Convergence behavior of the largest residual norm in the spectral window both with and without the k-means shift update.
    Overall residual convergence for the SISLICE method within the spectral window with (d) and without (c) k-means shift update.
    Converges in 2 SCF iterations both with and without k-means shift update.
  }
  \label{fig:silane_core_cluster_fixed}
\end{figure}

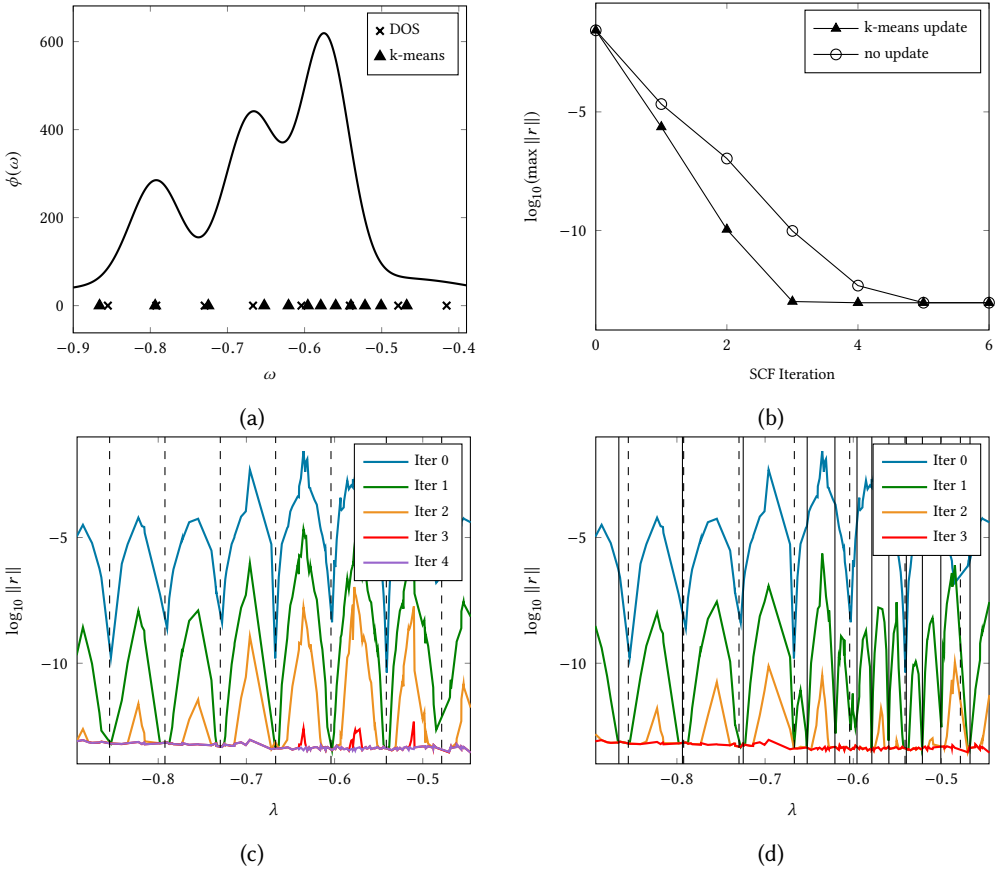
\begin{figure}[tbhp]

  \begin{subfigure}[b]{0.49\textwidth}
    \begin{tikzpicture}
\pgfmathsetlengthmacro\MajorTickLength{
  \pgfkeysvalueof{/pgfplots/major tick length} * 0.5
}
\begin{axis}[
  width=\textwidth,
  every axis plot/.append style={ thick},
  legend style={font=\tiny},
  legend cell align={left},
  every tick label/.append style={font=\tiny},
  major tick length=\MajorTickLength,
  xmin=-0.9, xmax=-0.39,
  ylabel={$\phi(\omega)$},
  xlabel={$\omega$},
  x label style={ font=\tiny },
  y label style={ font=\tiny }
]
\addplot[only marks, mark=x] table {
-0.85510383725 0
-0.79234521656 0
-0.72958659586 0
-0.66682797518 0
-0.60406935448 0
-0.54131073380 0
-0.47855211311 0
-0.41579349242 0
};
\addplot[only marks, mark=triangle*] table {
-0.86586008060 0
-0.79371613835 0
-0.72480726335 0
-0.65226613404 0
-0.62081221653 0
-0.59577895719 0
-0.57895438322 0
-0.55963087776 0
-0.53984862321 0
-0.52162384966 0
-0.50074073649 0
-0.46759324638 0
};
\legend{
DOS,
k-means
};
\addplot[mark={}] table [x=x, y=y, col sep=comma] {figs/SiOSi5_fixed_fermi_cluster_dos.dat};

\end{axis}
\end{tikzpicture}
    \caption{}
    \label{fig:silane_fermi_cluster_fixed_dos}
  \end{subfigure}~
  \begin{subfigure}[b]{0.49\textwidth}
    \begin{tikzpicture}
\pgfmathsetlengthmacro\MajorTickLength{
  \pgfkeysvalueof{/pgfplots/major tick length} * 0.5
}
\begin{axis}[
  width=\textwidth,
  legend style={font=\tiny},
  legend cell align={left},
  every tick label/.append style={font=\tiny},
  major tick length=\MajorTickLength,
  xmin=0, xmax=6,
  ylabel={$\log_{10}(\mathrm{max}\,\|r\|)$},
  xlabel={SCF Iteration},
  x label style={ font=\tiny },
  y label style={ font=\tiny }
]
\addplot[mark=triangle*,fill opacity=1] table [x=x, y=y3, col sep=comma] {figs/SiOSi5_fixed_waterfall.dat};
\addplot[mark=o,fill opacity=0] table [x=x, y=y4, col sep=comma] {figs/SiOSi5_fixed_waterfall.dat};
\legend{
k-means update, 
no update,
};
\end{axis}
\end{tikzpicture}
    \caption{}
    \label{fig:silane_fermi_cluster_fixed_waterfall}
  \end{subfigure}~\\

  \begin{subfigure}[b]{0.49\textwidth}
    \begin{tikzpicture}
\pgfmathsetlengthmacro\MajorTickLength{
  \pgfkeysvalueof{/pgfplots/major tick length} * 0.5
}
\begin{axis}[
  width=\textwidth,
  legend style={font=\tiny},
  legend cell align={left},
  every axis plot/.append style={ thick},
  every tick label/.append style={font=\tiny},
  major tick length=\MajorTickLength,
  xmin=-0.892, xmax=-0.446,
  ymin=-14 , ymax=-1,
  ylabel={$\log_{10}\|r\|$},
  xlabel={$\lambda$},
  x label style={ font=\tiny },
  y label style={ font=\tiny }
]
\addplot[mark={},color=cerulean] table [x=x, y=y1, col sep=comma] {figs/SiOSi5_fixed_fermi_cluster_conv_nokmeans.dat};
\addplot[mark={},color=applegreen] table [x=x, y=y2, col sep=comma] {figs/SiOSi5_fixed_fermi_cluster_conv_nokmeans.dat};
\addplot[mark={},color=carrot] table [x=x, y=y3, col sep=comma] {figs/SiOSi5_fixed_fermi_cluster_conv_nokmeans.dat};
\addplot[mark={},color=red] table [x=x, y=y4, col sep=comma] {figs/SiOSi5_fixed_fermi_cluster_conv_nokmeans.dat};
\addplot[mark={},color=amethyst] table [x=x, y=y5, col sep=comma] {figs/SiOSi5_fixed_fermi_cluster_conv_nokmeans.dat};
\draw[dashed] (-0.85510383725, -14) -- (-0.85510383725, -1);
\draw[dashed] (-0.79234521656, -14) -- (-0.79234521656, -1);
\draw[dashed] (-0.72958659586, -14) -- (-0.72958659586, -1);
\draw[dashed] (-0.66682797518, -14) -- (-0.66682797518, -1);
\draw[dashed] (-0.60406935448, -14) -- (-0.60406935448, -1);
\draw[dashed] (-0.54131073380, -14) -- (-0.54131073380, -1);
\draw[dashed] (-0.47855211311, -14) -- (-0.47855211311, -1);
\draw[dashed] (-0.41579349242, -14) -- (-0.41579349242, -1);
\legend{
Iter 0,
Iter 1,
Iter 2,
Iter 3,
Iter 4
};
\end{axis}
\end{tikzpicture}
    \caption{}
    \label{fig:silane_fermi_cluster_fixed_conv_nokmeans}
  \end{subfigure}~
  \begin{subfigure}[b]{0.49\textwidth}
    \begin{tikzpicture}
\pgfmathsetlengthmacro\MajorTickLength{
  \pgfkeysvalueof{/pgfplots/major tick length} * 0.5
}
\begin{axis}[
  width=\textwidth,
  legend style={font=\tiny},
  legend cell align={left},
  every axis plot/.append style={ thick},
  every tick label/.append style={font=\tiny},
  major tick length=\MajorTickLength,
  xmin=-0.892, xmax=-0.446,
  ymin=-14 , ymax=-1,
  ylabel={$\log_{10}\|r\|$},
  xlabel={$\lambda$},
  x label style={ font=\tiny },
  y label style={ font=\tiny }
]
\addplot[mark={},color=cerulean] table [x=x, y=y1, col sep=comma] {figs/SiOSi5_fixed_fermi_cluster_conv.dat};
\addplot[mark={},color=applegreen] table [x=x, y=y2, col sep=comma] {figs/SiOSi5_fixed_fermi_cluster_conv.dat};
\addplot[mark={},color=carrot] table [x=x, y=y3, col sep=comma] {figs/SiOSi5_fixed_fermi_cluster_conv.dat};
\addplot[mark={},color=red] table [x=x, y=y4, col sep=comma] {figs/SiOSi5_fixed_fermi_cluster_conv.dat};
\draw[dashed] (-0.85510383725, -14) -- (-0.85510383725, -1);
\draw[dashed] (-0.79234521656, -14) -- (-0.79234521656, -1);
\draw[dashed] (-0.72958659586, -14) -- (-0.72958659586, -1);
\draw[dashed] (-0.66682797518, -14) -- (-0.66682797518, -1);
\draw[dashed] (-0.60406935448, -14) -- (-0.60406935448, -1);
\draw[dashed] (-0.54131073380, -14) -- (-0.54131073380, -1);
\draw[dashed] (-0.47855211311, -14) -- (-0.47855211311, -1);
\draw[dashed] (-0.41579349242, -14) -- (-0.41579349242, -1);

\draw[] (-0.86586008060, -14) -- (-0.86586008060, -1);
\draw[] (-0.79371613835, -14) -- (-0.79371613835, -1);
\draw[] (-0.72480726335, -14) -- (-0.72480726335, -1);
\draw[] (-0.65226613404, -14) -- (-0.65226613404, -1);
\draw[] (-0.62081221653, -14) -- (-0.62081221653, -1);
\draw[] (-0.59577895719, -14) -- (-0.59577895719, -1);
\draw[] (-0.57895438322, -14) -- (-0.57895438322, -1);
\draw[] (-0.55963087776, -14) -- (-0.55963087776, -1);
\draw[] (-0.53984862321, -14) -- (-0.53984862321, -1);
\draw[] (-0.52162384966, -14) -- (-0.52162384966, -1);
\draw[] (-0.50074073649, -14) -- (-0.50074073649, -1);
\draw[] (-0.46759324638, -14) -- (-0.46759324638, -1);
\legend{
Iter 0,
Iter 1,
Iter 2,
Iter 3
};
\end{axis}
\end{tikzpicture}
    \caption{}
    \label{fig:silane_fermi_cluster_fixed_conv}
  \end{subfigure}

  \caption{Embedded cluster of Silane eigenvalues $[-0.9, -0.39]$ (141 eigenvalues) (a) Initial Lanczos DOS along with
    DOS shift placement and k-means update. (b) Convergence behavior of the largest residual norm in the spectral window both with and without the k-means shift update.
    Overall residual convergence for the SISLICE method within the spectral window with (d) and without (c) k-means shift update.
    Converges in 4 SCF iterations with k-means shift update and 6 iterations without shift update.
  }
  \label{fig:silane_fermi_cluster_fixed}
\end{figure}

\textbf{Silane}. For the case of Silane, SISLICE was applied to perform a full diagonalization using
100 shifts so that $\neig/\nsft \approx 11$.
We examine two representative spectral windows for this test case,
$C_1=[-20.59, -20.55]$ (\cref{fig:silane_core_cluster_fixed})
and $C_2=[-0.9, -0.39]$ (\cref{fig:silane_fermi_cluster_fixed}). The $C_1$ window represents a dense, isolated cluster of eigenvalues,
while the eigenvalues in $C_2$ are embedded in a dense region of eigenvalues.
Due to the different distribution characteristics of these two spectral windows, the convergence behavior of the
eigenpairs within these windows are different. However, because these windows are not treated separately
in the sense of the larger eigenvalue calculation, SCF iterations are performed until convergence is reached
across the spectrum.
Further, in this test case, DOS-based shift selection yielded 25 useless probes that were not well placed, i.e. probes
which did not produce any validated eigenvalues after the validation scheme outlined in \cref{sec:valid} was applied.
These probes were redistributed in the
subsequent iterations via the method outlined in \cref{sec:kmeans}.

Because eigenvalues in $C_1$ are well separated from the rest of the spectrum, the convergence of the subspace iteration is rapid.
Using the DOS based shift partitioning,
a single shift is placed just below $\lambda = -20.57$ to account for the 37 eigenvalues
in the immediate vicinity.
In the first SCF iteration, the eigenvalues
near the selected shift converged much more rapidly than those further away.
After the first SCF iteration, k-means eigenvalue clustering yielded
3 clusters of $\sim 12$ eigenvalues with centroids shown in \cref{fig:silane_core_cluster_fixed_dos}.
Convergence for this spectral window is achieved within 2 SCF iterations both with and without the k-means shift update, with all eigenvalues converging
at roughly the same rate notwithstanding their distance to the nearest shift. We can also see that for
the case of this isolated cluster, k-means clustering yielded no noticeable effects on residual convergence.

In contrast, the convergence of approximate eigenvalues in $C_2$ is less rapid due to the fact that there exist eigenvalues
both immediately below and above the eigenvalues in this spectral window. DOS based shift selection and spectrum partitioning
placed 8 evenly spaced shifts to account for the 141 eigenvalues within this window. After the first SCF iteration,
k-means clustering revealed a non-uniform distribution of eigenvalues within this window, yielding 12 clusters
of $\sim 11$ eigenvalues. Convergence rates for the eigenvalues in this window vary considerably based on their distance
to their nearest shift. Convergence across the entire spectral window is achieved within 4 SCF iterations with the k-means update
and 6 SCF iterations without the update. Thus, for this cluster, the k-means shift update yielded a discernible improvement in
the residual convergence of the approximate eigenpairs.

\textbf{Graphene}. For the case of Graphene, SISLICE was
applied to perform a partial diagonalization of the lowest $1000$ eigenvalues
using 100 shifts to obtain $\neig/\nsft \approx 10$. As can be seen in
\cref{fig:graphene_dos}, the eigenvalue distribution for Graphene is
approximately uniform. As such, the DOS based shift selection produced
uniformly distributed shifts along the entire spectral window, yielding no
useless probes.  We examine the eigenvalue interval $C = [ -1.4, -1.3 ]$ as a
representative example of the convergence behavior for this test case.

Within $C$, the Graphene test case admits 93 eigenvalues in a roughly uniform
distribution. As such, the DOS based shift partitioning places 8 evenly spaced
shifts in this spectral window so that $\neig/\nsft\approx 11$.
Within each spectral probe, convergence is more rapid near the shifts than
further away.
The k-means shift update simply migrates the shifts without any appreciable
changes to the shift spacing, i.e. the k-means result yields 8 clusters of
$\sim 11$ Ritz values with centroids of roughly equal spacing. Convergence for this
spectral window is achieved within 3 SCF iterations both with and without the k-means
shift update. We note that in the second SCF iteration, the convergence is slightly
worse with the k-means selected shift than the DOS selected shifts. However, as
both methods yield convergence in the same number of SCF iterations overall, we do
not believe this discrepancy to be problematic in practice.

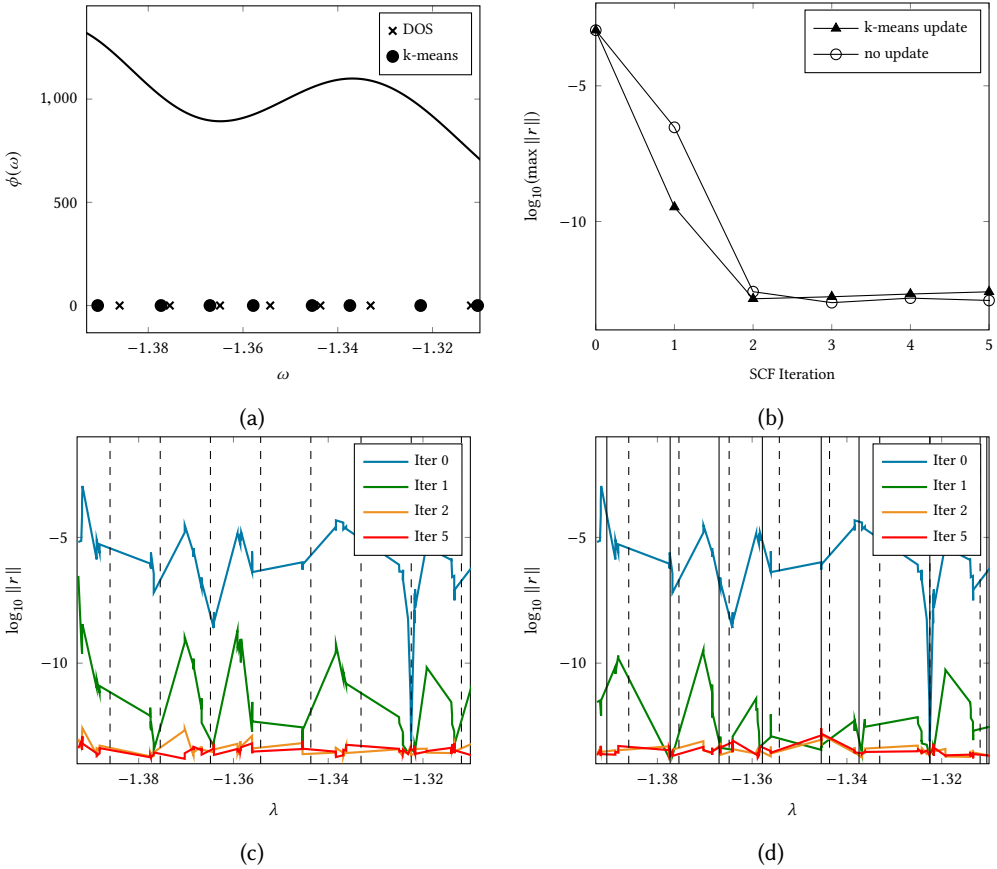
\begin{figure}[tbhp]
  \begin{subfigure}[b]{0.49\textwidth}
    \begin{tikzpicture}
\pgfmathsetlengthmacro\MajorTickLength{
  \pgfkeysvalueof{/pgfplots/major tick length} * 0.5
}
\begin{axis}[
  width=\textwidth,
  every axis plot/.append style={ thick},
  legend style={font=\tiny},
  legend cell align={left},
  every tick label/.append style={font=\tiny},
  major tick length=\MajorTickLength,
  xmin=-1.393, xmax=-1.31,
  ylabel={$\phi(\omega)$},
  xlabel={$\omega$},
  x label style={ font=\tiny },
  y label style={ font=\tiny }
]
\addplot[only marks, mark=x] table [
  x expr = {\thisrowno{0}},
  y expr = {0.0}
]
{
-1.396637923 
-1.386043864 
-1.375449804 
-1.364855745 
-1.354261685 
-1.343667626 
-1.333073567 
-1.322479507 
-1.311885448 
-1.301291388 
-1.290697329 
-1.280103269 
-1.269509210 
-1.258915151 
-1.248321091 
-1.237727032 
-1.227132973 
-1.216538913 
};

\addplot[only marks, mark=*] table [
  x expr = {\thisrowno{0}},
  y expr = {0.0}
]
{
-1.390670397
-1.377303788
-1.366982221
-1.357840859
-1.345417939
-1.337441872
-1.322501427
-1.310474319
-1.283132411
-1.260707958
-1.248119295
-1.230265100
-1.218216949
-1.202168590
};
\legend{
DOS,
k-means
};
\addplot[mark={}] table [x=x, y=y, col sep=comma] {figs/Graphene720_fixed_big_cluster_dos.dat};
\end{axis}
\end{tikzpicture}
    \caption{}
    \label{fig:graphene_big_cluster_fixed_dos}
  \end{subfigure}~
  \begin{subfigure}[b]{0.49\textwidth}
    \begin{tikzpicture}
\pgfmathsetlengthmacro\MajorTickLength{
  \pgfkeysvalueof{/pgfplots/major tick length} * 0.5
}
\begin{axis}[
  width=\textwidth,
  legend style={font=\tiny},
  legend cell align={left},
  every tick label/.append style={font=\tiny},
  major tick length=\MajorTickLength,
  xmin=0, xmax=5,
  ylabel={$\log_{10}(\mathrm{max}\,\|r\|)$},
  xlabel={SCF Iteration},
  x label style={ font=\tiny },
  y label style={ font=\tiny }
]
\addplot[mark=triangle*,fill opacity=1] table [x=x, y=y1, col sep=comma] {figs/Graphene720_fixed_waterfall.dat};
\addplot[mark=o,fill opacity=0] table [x=x, y=y2, col sep=comma] {figs/Graphene720_fixed_waterfall.dat};
\legend{
k-means update, 
no update,
};
\end{axis}
\end{tikzpicture}
    \caption{}
    \label{fig:graphene_big_cluster_fixed_waterfall}
  \end{subfigure}~\\

  \begin{subfigure}[b]{0.49\textwidth}
    \begin{tikzpicture}
\pgfmathsetlengthmacro\MajorTickLength{
  \pgfkeysvalueof{/pgfplots/major tick length} * 0.5
}
\begin{axis}[
  width=\textwidth,
  legend style={font=\tiny},
  legend cell align={left},
  every axis plot/.append style={ thick},
  every tick label/.append style={font=\tiny},
  major tick length=\MajorTickLength,
  xmin=-1.393, xmax=-1.31,
  ymin=-14, ymax=-1,
  xtick={-1.38,-1.36,-1.34,-1.32},
  ylabel={$\log_{10}\|r\|$},
  xlabel={$\lambda$},
  x label style={ font=\tiny },
  y label style={ font=\tiny }
]
\addplot[mark={},color=cerulean] table [x=x, y=y1, col sep=comma] {figs/Graphene720_fixed_big_cluster_conv_nokmeans.dat};
\addplot[mark={},color=applegreen] table [x=x, y=y2, col sep=comma] {figs/Graphene720_fixed_big_cluster_conv_nokmeans.dat};
\addplot[mark={},color=carrot] table [x=x, y=y3, col sep=comma] {figs/Graphene720_fixed_big_cluster_conv_nokmeans.dat};
\addplot[mark={},color=red] table [x=x, y=y4, col sep=comma] {figs/Graphene720_fixed_big_cluster_conv_nokmeans.dat};

\draw[dashed] (-1.396637923, -14) -- (-1.396637923 , -1);
\draw[dashed] (-1.386043864, -14) -- (-1.386043864 , -1);
\draw[dashed] (-1.375449804, -14) -- (-1.375449804 , -1);
\draw[dashed] (-1.364855745, -14) -- (-1.364855745 , -1);
\draw[dashed] (-1.354261685, -14) -- (-1.354261685 , -1);
\draw[dashed] (-1.343667626, -14) -- (-1.343667626 , -1);
\draw[dashed] (-1.333073567, -14) -- (-1.333073567 , -1);
\draw[dashed] (-1.322479507, -14) -- (-1.322479507 , -1);
\draw[dashed] (-1.311885448, -14) -- (-1.311885448 , -1);
\draw[dashed] (-1.301291388, -14) -- (-1.301291388 , -1);
\draw[dashed] (-1.290697329, -14) -- (-1.290697329 , -1);
\draw[dashed] (-1.280103269, -14) -- (-1.280103269 , -1);
\draw[dashed] (-1.269509210, -14) -- (-1.269509210 , -1);
\draw[dashed] (-1.258915151, -14) -- (-1.258915151 , -1);
\draw[dashed] (-1.248321091, -14) -- (-1.248321091 , -1);
\draw[dashed] (-1.237727032, -14) -- (-1.237727032 , -1);
\draw[dashed] (-1.227132973, -14) -- (-1.227132973 , -1);
\draw[dashed] (-1.216538913, -14) -- (-1.216538913 , -1);

\legend{
Iter 0,
Iter 1,
Iter 2,
Iter 5
};
\end{axis}
\end{tikzpicture}
    \caption{}
    \label{fig:graphene_big_cluster_fixed_conv_nokmeans}
  \end{subfigure}~
  \begin{subfigure}[b]{0.49\textwidth}
    \begin{tikzpicture}
\pgfmathsetlengthmacro\MajorTickLength{
  \pgfkeysvalueof{/pgfplots/major tick length} * 0.5
}
\begin{axis}[
  width=\textwidth,
  legend style={font=\tiny},
  legend cell align={left},
  every axis plot/.append style={ thick},
  every tick label/.append style={font=\tiny},
  major tick length=\MajorTickLength,
  xmin=-1.393, xmax=-1.31,
  ymin=-14, ymax=-1,
  xtick={-1.38,-1.36,-1.34,-1.32},
  ylabel={$\log_{10}\|r\|$},
  xlabel={$\lambda$},
  x label style={ font=\tiny },
  y label style={ font=\tiny }
]
\addplot[mark={},color=cerulean] table [x=x, y=y1, col sep=comma] {figs/Graphene720_fixed_big_cluster_conv.dat};
\addplot[mark={},color=applegreen] table [x=x, y=y2, col sep=comma] {figs/Graphene720_fixed_big_cluster_conv.dat};
\addplot[mark={},color=carrot] table [x=x, y=y3, col sep=comma] {figs/Graphene720_fixed_big_cluster_conv.dat};
\addplot[mark={},color=red] table [x=x, y=y4, col sep=comma] {figs/Graphene720_fixed_big_cluster_conv.dat};

\draw[dashed] (-1.396637923, -14) -- (-1.396637923 , -1);
\draw[dashed] (-1.386043864, -14) -- (-1.386043864 , -1);
\draw[dashed] (-1.375449804, -14) -- (-1.375449804 , -1);
\draw[dashed] (-1.364855745, -14) -- (-1.364855745 , -1);
\draw[dashed] (-1.354261685, -14) -- (-1.354261685 , -1);
\draw[dashed] (-1.343667626, -14) -- (-1.343667626 , -1);
\draw[dashed] (-1.333073567, -14) -- (-1.333073567 , -1);
\draw[dashed] (-1.322479507, -14) -- (-1.322479507 , -1);
\draw[dashed] (-1.311885448, -14) -- (-1.311885448 , -1);
\draw[dashed] (-1.301291388, -14) -- (-1.301291388 , -1);
\draw[dashed] (-1.290697329, -14) -- (-1.290697329 , -1);
\draw[dashed] (-1.280103269, -14) -- (-1.280103269 , -1);
\draw[dashed] (-1.269509210, -14) -- (-1.269509210 , -1);
\draw[dashed] (-1.258915151, -14) -- (-1.258915151 , -1);
\draw[dashed] (-1.248321091, -14) -- (-1.248321091 , -1);
\draw[dashed] (-1.237727032, -14) -- (-1.237727032 , -1);
\draw[dashed] (-1.227132973, -14) -- (-1.227132973 , -1);
\draw[dashed] (-1.216538913, -14) -- (-1.216538913 , -1);

\draw[] (-1.390670397, -14) -- (-1.390670397, -1);
\draw[] (-1.377303788, -14) -- (-1.377303788, -1);
\draw[] (-1.366982221, -14) -- (-1.366982221, -1);
\draw[] (-1.357840859, -14) -- (-1.357840859, -1);
\draw[] (-1.345417939, -14) -- (-1.345417939, -1);
\draw[] (-1.337441872, -14) -- (-1.337441872, -1);
\draw[] (-1.322501427, -14) -- (-1.322501427, -1);
\draw[] (-1.310474319, -14) -- (-1.310474319, -1);
\draw[] (-1.283132411, -14) -- (-1.283132411, -1);
\draw[] (-1.260707958, -14) -- (-1.260707958, -1);
\draw[] (-1.248119295, -14) -- (-1.248119295, -1);
\draw[] (-1.230265100, -14) -- (-1.230265100, -1);
\draw[] (-1.218216949, -14) -- (-1.218216949, -1);
\draw[] (-1.202168590, -14) -- (-1.202168590, -1);

\legend{
Iter 0,
Iter 1,
Iter 2,
Iter 5
};
\end{axis}
\end{tikzpicture}
    \caption{}
    \label{fig:graphene_big_cluster_fixed_conv}
  \end{subfigure}

  \caption{Graphene eigenvalue cluster $[-1.4,-1.3]$ (93 eigenvalues). (a) Initial Lanczos DOS along with
    DOS shift placement and k-means update. (b) Convergence behavior of the largest residual norm in the spectral window both with and without the k-means shift update.
    Overall residual convergence for the SISLICE method within the spectral window with (d) and without (c) k-means shift update.
    Converges in 3 SCF iterations both with and without k-means shift update.
  }
  \label{fig:graphene_big_cluster_fixed}
\end{figure}

\subsection{Shift Selection for a Converging Matrix Pencil Sequence}

In this section, we examine how our shift selection strategy enables
the SISLICE method to efficiently compute eigenpairs of pregenerated, convergent sequences of
matrix pencils obtained from a true SCF procedure. Such a test allows us to
gauge
the ability of the SISLICE method to solve true SCF eigenvalue problems. To
determine the efficacy of the shift selection and migration strategy, we examine both the
convergence
of the residuals produced by the SISLICE method and the change of the true
eigenvalues throughout the SCF procedure itself. The latter is possible because
these matrices are pregenerated, thus we have access to the exact eigenvalues
of these matrices as a reference to compare the convergence of the SISLICE
method.
Further, as was examined in the previous section, we perform analogous
comparisons of the SISLICE method both with and without k-means updates to the
spectral shifts throughout the SCF procedure.

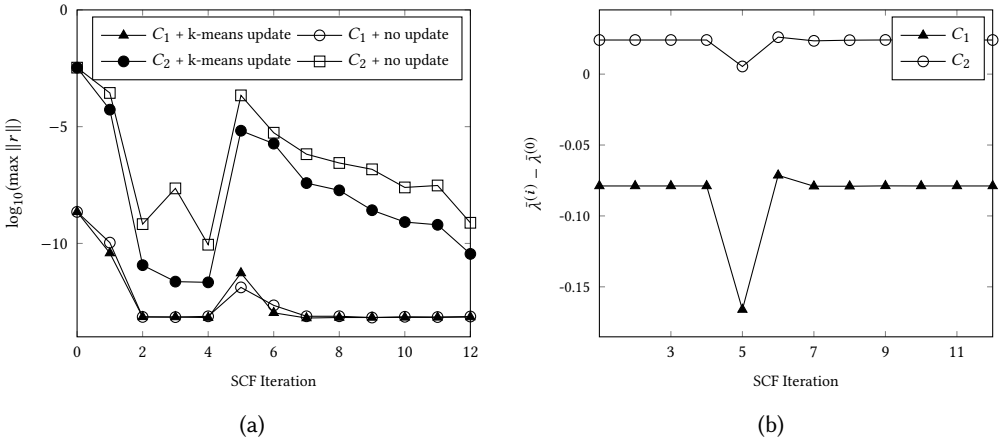
\begin{figure}[t]
  \begin{subfigure}[b]{0.49\textwidth}
    \begin{tikzpicture}
\pgfmathsetlengthmacro\MajorTickLength{
  \pgfkeysvalueof{/pgfplots/major tick length} * 0.5
}
\begin{axis}[
  width=\textwidth,
  legend style={font=\tiny},
  legend cell align={left},
  legend columns=2,
  every tick label/.append style={font=\tiny},
  major tick length=\MajorTickLength,
  xmin=0, xmax=12,
  ymin=-14, ymax=0,
  ylabel={$\log_{10}(\mathrm{max}\,\|r\|)$},
  xlabel={SCF Iteration},
  x label style={ font=\tiny },
  y label style={ font=\tiny }
]
\addplot[mark=triangle*,fill opacity=1] table [x=x, y=y1, col sep=comma] {figs/SiOSi5_dynamic_waterfall.dat};
\addplot[mark=o,fill opacity=0] table [x=x, y=y3, col sep=comma] {figs/SiOSi5_dynamic_waterfall.dat};
\addplot[mark=*,fill opacity=1] table [x=x, y=y2, col sep=comma] {figs/SiOSi5_dynamic_waterfall.dat};
\addplot[mark=square*,fill opacity=0] table [x=x, y=y4, col sep=comma] {figs/SiOSi5_dynamic_waterfall.dat};
\legend{
$C_1$ + k-means update, 
$C_1$ + no update,
$C_2$ + k-means update, 
$C_2$ + no update
};
\end{axis}
\end{tikzpicture}
    \caption{}
    \label{fig:silane_dynamic_waterfall}
  \end{subfigure}~
  \begin{subfigure}[b]{0.49\textwidth}
    \begin{tikzpicture}
\pgfmathsetlengthmacro\MajorTickLength{
  \pgfkeysvalueof{/pgfplots/major tick length} * 0.5
}
\begin{axis}[
  width=\textwidth,
  legend style={font=\tiny},
  legend cell align={left},
  every tick label/.append style={font=\tiny},
  major tick length=\MajorTickLength,
  xmin=1, xmax=12,
  ylabel={$\bar\lambda^{(i)} - \bar\lambda^{(0)}$},
  xlabel={SCF Iteration},
  ytick={0,-0.05,-0.10,-0.15},
  yticklabels={0,-0.05,-0.10,-0.15},
  xtick={3,5,7,9,11},
  x label style={ font=\tiny },
  y label style={ font=\tiny }
]
\addplot[mark=triangle*,fill opacity=1] table [x=x, y=y1, col sep=comma] {figs/SiOSi5_dynamic_true_eig.dat};
\addplot[mark=o,fill opacity=0] table [x=x, y=y2, col sep=comma] {figs/SiOSi5_dynamic_true_eig.dat};
\legend{
$C_1$, 
$C_2$ 
};
\end{axis}
\end{tikzpicture}
    \caption{}
    \label{fig:silane_dynamic_true_eig}
  \end{subfigure}
  \caption{Convergence of the SISLICE method applied to the Silane SCF procedure for
    two representative spectral windows. (a) Comparison of the convergence behavior of the largest residual
    norm in the respective spectral windows both with and without k-means shift updates throughout the
    SCF procedure. (b) The change in average eigenvalue from the initial average of the two spectral
    windows respectively.
  }
  \label{fig:silane_dynamic}
\end{figure}

\textbf{Silane}. The Silane SCF procedure converged within 13 iterations in the
NWChemEx software package.
It is the nature of
this particular test case (and is typical of all-electron density functional
theory calculation) that the spectrum is separated into well defined clusters
throughout the entire SCF procedure. For this reason, we are able to examine
the same eigenvalue clusters as discussed in the previous section for this test
case. The convergence behavior of the SISLICE method applied to this test case
is given in \cref{fig:silane_dynamic}.

Much like the results presented in \cref{sec:fixed}, we see
a significant difference in the convergence behavior between the two clusters.
Due to the isolated nature of $C_1$, convergence of the
subspace iteration is rapid. Despite changes in
the eigenvalues resulting from the changes in the matrix pencil in early
SCF iterations, the shift selection strategy we developed is able to
track this change, and construct and move spectral probes to
obtain eigenvalues within $C_1$ at convergence.
The convergence for $C_2$ is much less rapid, as was
also the case in the previous experiments due to the lack of a large
separation between eigenvalues within $C_2$ and the rest of the spectrum.
Furthermore, the change in the eigenvectors in $C_1$ is
much less than those in $C_2$, thus they provide excellent initial guesses
for subsequent SCF iterations. The eigenvectors in $C_2$ undergo a much more drastic
change, but it can be seen in \cref{fig:silane_dynamic_waterfall} that
this change becomes less as the SCF converges.


Note that the convergence of the SCF for the Silane test case is not smooth;
there is a large change in the average eigenvalue for the two examined clusters at
the fifth SCF iteration. This is not an uncommon feature in the SCF procedure
for density functional calculations. There is an analogous change in the residual
norms for the SISLICE method at the same SCF iteration. The reason for this is
two-fold. In the case where k-means clustering is used to migrate the shifts
between SCF iterations, the fact that the clustering is performed using the validated
eigenpairs from the \emph{previous} SCF iteration yielded a non-optimal placement for the 5th iteration.
However, because the change in residual norms is present also for the experiment
without k-means shift updates, the shift migration is not the only reason for this change.
The large change in average eigenvalue for this SCF iteration is also accompanied by
a change in corresponding eigenvectors within these spectral windows. Thus, the validated
eigenvectors from the previous SCF iteration are also not optimal choices for the initial
guess to seed the subspace sequence at this iteration. Remark that the increase in residual
norm is in fact less for the $C_2$ cluster with k-means shift updates, indicating that the
shift migration strategy is beneficial for this cluster even when the shifts are placed
non-optimally.

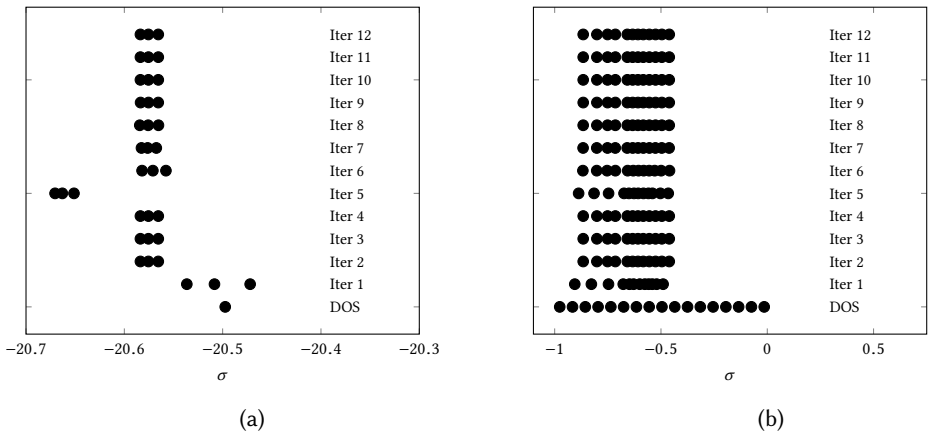
\begin{figure}[t]
  \centering
  \begin{subfigure}[b]{0.49\textwidth}
    \begin{tikzpicture}
\pgfmathsetlengthmacro\MajorTickLength{
  \pgfkeysvalueof{/pgfplots/major tick length} * 0.5
}
\begin{axis}[
  width=\textwidth,
  legend style={font=\tiny},
  legend cell align={left},
  every tick label/.append style={font=\tiny},
  major tick length=\MajorTickLength,
  xmin=-20.7, xmax=-20.3,
  ytick={},
  yticklabels={},
  xlabel={$\sigma$},
  x label style={ font=\tiny },
  y label style={ font=\tiny }
]
\addplot[only marks] table {
  -20.497386192 0
};
\addplot[only marks] table {
  -20.536444924 1
  -20.508324550 1
  -20.472035967 1
};
\addplot[only marks] table {
  -2.058365831700000115e+01 2
  -2.057550221800000045e+01 2
  -2.056546808099999879e+01 2
};
\addplot[only marks] table {
  -2.058365831700000115e+01 3
  -2.057550221800000045e+01 3
  -2.056546808099999879e+01 3
};
\addplot[only marks] table {
  -2.058365831700000115e+01 4
  -2.057550221800000045e+01 4
  -2.056546808099999879e+01 4
};
\addplot[only marks] table {
  -2.067046420499999826e+01 5
  -2.066303570099999831e+01 5
  -2.065121625700000152e+01 5
};
\addplot[only marks] table {
  -2.058198787099999905e+01 6
  -2.057096218099999874e+01 6
  -2.055769260699999990e+01 6
};
\addplot[only marks] table {
  -2.058258725399999989e+01 7
  -2.057630460100000036e+01 7
  -2.056748370000000037e+01 7
};
\addplot[only marks] table {
  -2.058430869900000104e+01 8
  -2.057576060699999942e+01 8
  -2.056530121600000172e+01 8
};
\addplot[only marks] table {
  -2.058348743399999847e+01 9
  -2.057539430899999999e+01 9
  -2.056549983400000059e+01 9
};
\addplot[only marks] table {
  -2.058365831700000115e+01 10 
  -2.057550221800000045e+01 10
  -2.056546808099999879e+01 10
};
\addplot[only marks] table {
  -2.058365831700000115e+01 11 
  -2.057550221800000045e+01 11
  -2.056546808099999879e+01 11
};
\addplot[only marks] table {
  -2.058365831700000115e+01 12 
  -2.057550221800000045e+01 12
  -2.056546808099999879e+01 12
};

\node[right,text=black] at (-20.4,0 )   {\tiny DOS};
\node[right,text=black] at (-20.4,1 )   {\tiny Iter 1};
\node[right,text=black] at (-20.4,2 )   {\tiny Iter 2 };
\node[right,text=black] at (-20.4,3 )   {\tiny Iter 3 };
\node[right,text=black] at (-20.4,4 )   {\tiny Iter 4 };
\node[right,text=black] at (-20.4,5 )   {\tiny Iter 5 };
\node[right,text=black] at (-20.4,6 )   {\tiny Iter 6 };
\node[right,text=black] at (-20.4,7 )   {\tiny Iter 7 };
\node[right,text=black] at (-20.4,8 )   {\tiny Iter 8 };
\node[right,text=black] at (-20.4,9 )   {\tiny Iter 9 };
\node[right,text=black] at (-20.4,10)   {\tiny Iter 10};
\node[right,text=black] at (-20.4,11)   {\tiny Iter 11};
\node[right,text=black] at (-20.4,12)   {\tiny Iter 12};
\end{axis}
\end{tikzpicture}
    \caption{}
    \label{fig:silane_dynamic_core_shift_migration}
  \end{subfigure}~
  \begin{subfigure}[b]{0.49\textwidth}
    \begin{tikzpicture}
\pgfmathsetlengthmacro\MajorTickLength{
  \pgfkeysvalueof{/pgfplots/major tick length} * 0.5
}
\begin{axis}[
  width=\textwidth,
  legend style={font=\tiny},
  legend cell align={left},
  every tick label/.append style={font=\tiny},
  major tick length=\MajorTickLength,
  xmin=-1.1, xmax=0.75,
  ytick={},
  yticklabels={},
  xlabel={$\sigma$},
  x label style={ font=\tiny },
  y label style={ font=\tiny }
]
\addplot[only marks] table {
-9.762498282099999480e-01 0
-9.161361918499999479e-01 0
-8.560225554899999478e-01 0
-7.959089191199999469e-01 0
-7.357952827599999468e-01 0
-6.756816463999999467e-01 0
-6.155680100299999458e-01 0
-5.554543736699999457e-01 0
-4.953407373000000002e-01 0
-4.352271009400000001e-01 0
-3.751134645800000000e-01 0
-3.149998282099999991e-01 0
-2.548861918499999990e-01 0
-1.947725554899999989e-01 0
-1.346589191199999980e-01 0
-7.454528275900000611e-02 0
-1.443164639500000082e-02 0
};
\addplot[only marks] table {
-9.054821020000000109e-01 1
-8.274663201299999882e-01 1
-7.468801171900000258e-01 1
-6.775953813399999692e-01 1
-6.484166188700000166e-01 1
-6.285542776600000492e-01 1
-5.999798269699999631e-01 1
-5.760149528600000401e-01 1
-5.576894527000000323e-01 1
-5.407773402600000168e-01 1
-5.201021208699999621e-01 1
-4.899476860900000008e-01 1
};
\addplot[only marks] table {
-8.658618214500000132e-01 2
-8.015865222800000467e-01 2
-7.512935881200000221e-01 2
-7.145628528399999846e-01 2
-6.583209120399999836e-01 2
-6.332671045300000134e-01 2
-6.082697319800000058e-01 2
-5.834276852799999968e-01 2
-5.549489356900000336e-01 2
-5.255825128700000182e-01 2
-4.957987247599999758e-01 2
-4.610757286900000196e-01 2
};
\addplot[only marks] table {
-8.658618214500000132e-01 3
-8.015865222800000467e-01 3
-7.512935881200000221e-01 3
-7.145628528399999846e-01 3
-6.583209120399999836e-01 3
-6.332671045300000134e-01 3
-6.082697319800000058e-01 3
-5.834276852799999968e-01 3
-5.549489356900000336e-01 3
-5.255825128700000182e-01 3
-4.957987247599999758e-01 3
-4.610757286900000196e-01 3
};
\addplot[only marks] table {
-8.658618214500000132e-01 4
-8.015865222800000467e-01 4
-7.512935881200000221e-01 4
-7.145628528399999846e-01 4
-6.583209120399999836e-01 4
-6.332671045300000134e-01 4
-6.082697319800000058e-01 4
-5.834276852799999968e-01 4
-5.549489356900000336e-01 4
-5.255825128700000182e-01 4
-4.957987247599999758e-01 4
-4.610757286900000196e-01 4
};
\addplot[only marks] table {
-8.873624993500000002e-01 5
-8.148891803400000189e-01 5
-7.464531167500000342e-01 5
-6.738820916200000122e-01 5
-6.508284386499999519e-01 5
-6.269915265400000104e-01 5
-6.069880211399999448e-01 5
-5.825184380600000322e-01 5
-5.589784387399999721e-01 5
-5.400790168300000360e-01 5
-5.029068893600000534e-01 5
-4.657645966499999735e-01 5
};
\addplot[only marks] table {
-8.674014236899999508e-01 6
-8.022804458100000202e-01 6
-7.534543485300000398e-01 6
-7.144346315399999758e-01 6
-6.568111178100000114e-01 6
-6.305512895399999485e-01 6
-6.056266712300000243e-01 6
-5.813544541699999790e-01 6
-5.557197053500000372e-01 6
-5.295188499699999829e-01 6
-4.993827531100000217e-01 6
-4.599174938400000245e-01 6
};

\addplot[only marks] table {
-8.664058195899999593e-01 7
-8.020606780500000310e-01 7
-7.516731265899999626e-01 7
-7.149835674900000271e-01 7
-6.590340062499999707e-01 7
-6.329409585999999921e-01 7
-6.079493471700000429e-01 7
-5.830715343400000306e-01 7
-5.547204631900000127e-01 7
-5.262880716100000189e-01 7
-4.962956198199999780e-01 7
-4.615564951100000202e-01 7
};
\addplot[only marks] table {
-8.659620371599999755e-01 8
-8.016628641999999916e-01 8
-7.513803231500000290e-01 8
-7.145896422099999867e-01 8
-6.584572085200000036e-01 8
-6.334189789599999898e-01 8
-6.084473977400000066e-01 8
-5.835491116000000034e-01 8
-5.550886800000000010e-01 8
-5.257414949499999990e-01 8
-4.959246345600000216e-01 8
-4.611430238099999968e-01 8
};
\addplot[only marks] table {
-8.657839191199999673e-01 9
-8.015117488699999848e-01 9
-7.512134566400000546e-01 9
-7.144880477999999702e-01 9
-6.582876092399999601e-01 9
-6.332266771800000305e-01 9
-6.082248841999999822e-01 9
-5.834029368299999962e-01 9
-5.549177658799999691e-01 9
-5.255545245099999541e-01 9
-4.957683557099999971e-01 9
-4.610720178600000163e-01 9
};
\addplot[only marks] table {
-8.658618214500000132e-01 10
-8.015865222800000467e-01 10
-7.512935881200000221e-01 10
-7.145628528399999846e-01 10
-6.583209120399999836e-01 10
-6.332671045300000134e-01 10
-6.082697319800000058e-01 10
-5.834276852799999968e-01 10
-5.549489356900000336e-01 10
-5.255825128700000182e-01 10
-4.957987247599999758e-01 10
-4.610757286900000196e-01 10
};
\addplot[only marks] table {
-8.658618214500000132e-01 11
-8.015865222800000467e-01 11
-7.512935881200000221e-01 11
-7.145628528399999846e-01 11
-6.583209120399999836e-01 11
-6.332671045300000134e-01 11
-6.082697319800000058e-01 11
-5.834276852799999968e-01 11
-5.549489356900000336e-01 11
-5.255825128700000182e-01 11
-4.957987247599999758e-01 11
-4.610757286900000196e-01 11
};
\addplot[only marks] table {
-8.658618214500000132e-01 12
-8.015865222800000467e-01 12
-7.512935881200000221e-01 12
-7.145628528399999846e-01 12
-6.583209120399999836e-01 12
-6.332671045300000134e-01 12
-6.082697319800000058e-01 12
-5.834276852799999968e-01 12
-5.549489356900000336e-01 12
-5.255825128700000182e-01 12
-4.957987247599999758e-01 12
-4.610757286900000196e-01 12
};

\node[right,text=black] at (0.25,0 )   {\tiny DOS};
\node[right,text=black] at (0.25,1 )   {\tiny Iter 1};
\node[right,text=black] at (0.25,2 )   {\tiny Iter 2 };
\node[right,text=black] at (0.25,3 )   {\tiny Iter 3 };
\node[right,text=black] at (0.25,4 )   {\tiny Iter 4 };
\node[right,text=black] at (0.25,5 )   {\tiny Iter 5 };
\node[right,text=black] at (0.25,6 )   {\tiny Iter 6 };
\node[right,text=black] at (0.25,7 )   {\tiny Iter 7 };
\node[right,text=black] at (0.25,8 )   {\tiny Iter 8 };
\node[right,text=black] at (0.25,9 )   {\tiny Iter 9 };
\node[right,text=black] at (0.25,10)   {\tiny Iter 10};
\node[right,text=black] at (0.25,11)   {\tiny Iter 11};
\node[right,text=black] at (0.25,12)   {\tiny Iter 12};
\end{axis}
\end{tikzpicture}
    \caption{}
    \label{fig:silane_dynamic_fermi_shift_migration}
  \end{subfigure}
  \caption{Shift migration for the SISLICE method in the Silane SCF procedure for the $C_1$ (a) and $C_2$ (b)
    eigenvalue clusters.}
  \label{fig:silane_shift_migration}
\end{figure}

Due to the fact that Silane admits well-defined (and trackable) clusters in its
spectrum, we are also able to examine the shift migration within these clusters
in \cref{fig:silane_shift_migration}. The largest change in shift placement
occurred between the first and second SCF iterations, the former of which was
produced by the DOS strategy. Because the characteristic of these clusters is largely
unchanging throughout the SCF procedure, we are able to see that the k-means
shift update remains visually unchanging with the exception of the fifth SCF iteration.
Remark that the k-means shift selection strategy was able to track the change in
eigenvalues in this iteration and subsequently recover to a reasonably static set
of shifts in the following SCF iterations.

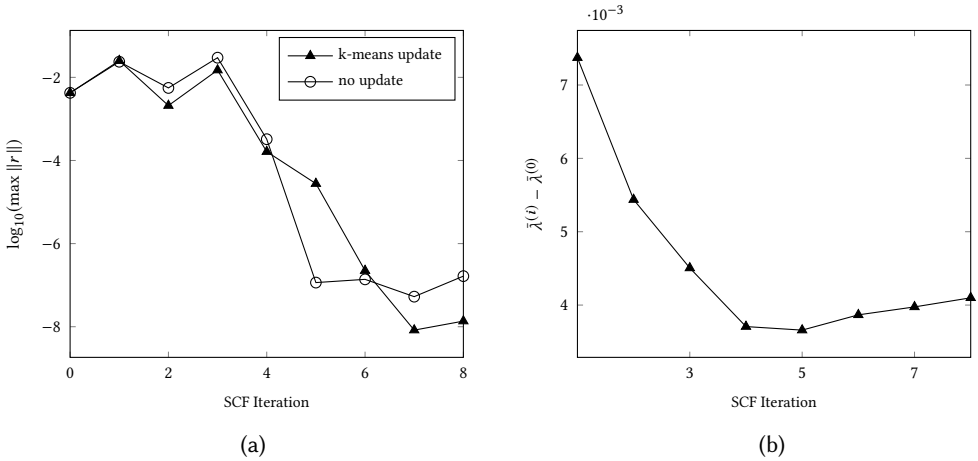
\begin{figure}[t]
  \begin{subfigure}[b]{0.49\textwidth}
    \begin{tikzpicture}
\pgfmathsetlengthmacro\MajorTickLength{
  \pgfkeysvalueof{/pgfplots/major tick length} * 0.5
}
\begin{axis}[
  width=\textwidth,
  legend style={font=\tiny},
  legend cell align={left},
  every tick label/.append style={font=\tiny},
  major tick length=\MajorTickLength,
  xmin=0, xmax=8,
  ylabel={$\log_{10}(\mathrm{max}\,\|r\|)$},
  xlabel={SCF Iteration},
  x label style={ font=\tiny },
  y label style={ font=\tiny }
]
\addplot[mark=triangle*,fill opacity=1] table [x=x, y=y1, col sep=comma] {figs/Graphene720_dynamic_waterfall.dat};
\addplot[mark=o,fill opacity=0] table [x=x, y=y2, col sep=comma] {figs/Graphene720_dynamic_waterfall.dat};
\legend{
k-means update, 
no update,
};
\end{axis}
\end{tikzpicture}
    \caption{}
    \label{fig:graphene_dynamic_waterfall}
  \end{subfigure}~
  \begin{subfigure}[b]{0.49\textwidth}
    \begin{tikzpicture}
\pgfmathsetlengthmacro\MajorTickLength{
  \pgfkeysvalueof{/pgfplots/major tick length} * 0.5
}
\begin{axis}[
  width=\textwidth,
  legend style={font=\tiny},
  legend cell align={left},
  every tick label/.append style={font=\tiny},
  major tick length=\MajorTickLength,
  xmin=1, xmax=8,
  ylabel={$\bar\lambda^{(i)} - \bar\lambda^{(0)}$},
  xlabel={SCF Iteration},
  xtick={3,5,7},
  x label style={ font=\tiny },
  y label style={ font=\tiny }
]
\addplot[mark=triangle*,fill opacity=1] table [x=x, y=y, col sep=comma] {figs/Graphene720_dynamic_true_eig.dat};
\end{axis}
\end{tikzpicture}
    \caption{}
    \label{fig:graphene_dynamic_true_eig}
  \end{subfigure}
  \caption{Convergence of the SISLICE method applied to the Graphene SCF procedure for the lowest 1000 eigenvalues.
    (a) Comparison of the convergence behavior of the largest residual
    norm in both with and without k-means shift updates throughout the
    SCF procedure. (b) The change in average eigenvalue from the initial average of the lowest 1000 eigenvalues.
  }
  \label{fig:graphene_dynamic}
\end{figure}

\textbf{Graphene}. The Graphene SCF procedure converged within 10 iterations in
the SIESTA software package. Convergence results for the SISLICE method applied
to this problem may be found in \cref{fig:graphene_dynamic}. Unlike the Silane
test case, the homogeneity of the eigenvectors for the Graphene test case makes
tracking eigenvalue clusters throughout the SCF procedure impractical. The
eigenpairs in a particular spectral interval at one SCF iteration are not likely
to be of the same character in the subsequent iterations until convergence is
reached. As such, we examine the convergence globally across all of the 1000
eigenpairs obtained desired for this test case.

Unlike the Silane test case, the SCF convergence for Graphene is smooth.
This smooth SCF convergence is mirrored in the monotonic convergence of the
SISLICE method as the SCF approaches convergence. When the SCF procedure yielded
large changes in the underlying spectrum, i.e. the first 3 iterations, the error produced by the SISLICE method
was larger as the bases from the previous SCF iteration were not as good of an initial
guess as they were in the later iterations. After the fourth SCF iteration, the spectrum is only
undergoing small changes and the SISLICE method exhibits rapid convergence. As was the case for the
previous numerical experiments with Graphene, no discernible difference
between DOS and k-means shifts is exhibited. For example, at the fifth SCF iteration,
the DOS shifts produced more accurate results whereas at the seventh, the k-means
shifts produced more accurate results.

Due to the fact that the SCF underwent large spectral changes in the early SCF iterations,
the extent to which the shifts were able to be usefully updated using the spectrum of the previous
matrix pencils was limited. As such, shifts needed to be inserted per the prescription in
\cref{sec:probe_insert}. In the following subsection, we examine an example of this insertion
for the Graphene test case.

\subsection{Missing Eigenvalues and Probe Insertion}
\label{sec:misev_results}

\begin{figure}[t]
\begin{subfigure}[b]{0.49\textwidth}
  \begin{tikzpicture}
\begin{axis}[ 
    every tick label/.append style={font=\tiny},
    xmin=-0.75, xmax=-0.56,
    ymin=-1.0, ymax=2.0,
    xlabel={\tiny $\lambda$},
    legend columns=2,
    ytick=\empty,
    width=\textwidth
]

\addplot[only marks, mark=*] table [
  x expr = {\thisrowno{0}},
  y expr = {0.5}
] {
-7.3883091259e-01
};

\addplot[only marks, mark=triangle*] table [
  x expr = {\thisrowno{0}},
  y expr = {0.5}
] {
-7.1855153121e-01
-6.9827214983e-01
-6.5771338707e-01
-6.1715462431e-01
-5.7659586155e-01
};

\addplot[only marks, mark=x, color=black] table [
  x expr = {\thisrowno{0}},
  y expr = {1.0}
] {
   -7.3524681565e-01
   -7.3512881238e-01
   -7.3512865622e-01
   -7.3512848248e-01
   -7.3512722858e-01
   -7.3424106578e-01
   -7.3423570618e-01
   -7.3421878204e-01
   -7.3416304022e-01
   -7.3323107385e-01
   -7.3322978601e-01
   -7.3322707496e-01
   -7.3322605605e-01
   -7.3084862465e-01
   -7.2936798437e-01
   -7.2547202543e-01
   -7.2335441445e-01
   -7.2235684112e-01
   -7.2208364583e-01
   -7.2180190692e-01
   -7.2156137603e-01
   -7.2061545507e-01
   -7.2003761701e-01
   -7.1874177753e-01
   -7.1598346320e-01
   -7.1349929309e-01
   -7.1081465466e-01
   -7.0917678622e-01
   -7.0666310270e-01
   -7.0558694679e-01
   -7.0054639924e-01
   -6.9910892410e-01
   -6.9557213002e-01
   -6.9191494058e-01
   -6.9080550214e-01
   -6.8623025572e-01
   -6.8408662492e-01
   -6.8093709335e-01
   -6.7997779649e-01
   -6.7961389149e-01
   -6.7567785163e-01
   -6.6982783643e-01
   -6.6381869022e-01
   -6.6327099825e-01
   -6.5900319901e-01
   -6.5595645902e-01
   -6.5526940884e-01
   -6.5449015369e-01
   -6.5401910828e-01
   -6.5333738239e-01
   -6.5200148852e-01
   -6.5172616229e-01
   -6.4534761530e-01
   -6.4364916000e-01
   -6.4357368673e-01
   -6.4357364059e-01
   -6.4355122834e-01
   -6.3992878831e-01
   -6.3844387837e-01
   -6.3789231343e-01
   -6.3766847962e-01
   -6.3487346099e-01
   -6.3340468976e-01
   -6.3246949239e-01
   -6.3199240134e-01
   -6.3186947416e-01
   -6.3162105715e-01
   -6.2965187214e-01
   -5.8932343587e-01
   -5.8922139252e-01
   -5.8120726702e-01
   -5.7818452458e-01
   -5.7694211019e-01
   -5.7669586155e-01
};

\addplot[only marks,mark=x,color=red] table [
  x expr = {\thisrowno{0}},
  y expr = {0.0}
] {
-7.3512881238e-01
-7.3512865622e-01
-7.3512848248e-01
-7.3512722858e-01
-7.3424106578e-01
-7.3423570618e-01
-7.3421878204e-01
-7.3416304022e-01
-7.3323107385e-01
-7.3322978601e-01
-7.3322707496e-01
-7.3322605605e-01
-7.2674984216e-01
-7.2673834840e-01
-7.2664222536e-01
-7.2651395899e-01
-7.2498184894e-01
-7.2496869995e-01
-7.2262230464e-01
-7.2258808300e-01
-7.1834565559e-01
-7.1834415901e-01
-7.1782426377e-01
-7.1780434571e-01
-7.1778709938e-01
-7.1776785051e-01
-7.1761637187e-01
-7.1761628922e-01
-7.1759070700e-01
-7.1758237129e-01
-7.1758230049e-01
-7.1756200854e-01
-7.1751744924e-01
-7.1748785417e-01
-7.1666568324e-01
-7.1666563055e-01
-7.1562437944e-01
-7.1562436320e-01
-7.1559757991e-01
-7.1559752835e-01
-7.1413344335e-01
-7.1413328474e-01
-7.1399314334e-01
-7.1399312496e-01
-7.1399174681e-01
-7.1399153910e-01
-7.1397855184e-01
-7.1397849650e-01
-7.1229245808e-01
-7.1228734513e-01
-7.1216931879e-01
-7.1216603865e-01
-7.1201312715e-01
-7.1200850999e-01
-7.1200765862e-01
-7.1200698770e-01
-7.1168030929e-01
-7.1167270847e-01
-7.1167055717e-01
-7.1166840801e-01
-7.1006600472e-01
-7.1000599894e-01
-7.0997349570e-01
-7.0996842172e-01
-7.0649762288e-01
-7.0647559386e-01
-7.0635445958e-01
-7.0633728837e-01
-7.0608166585e-01
-7.0606309511e-01
-7.0605261503e-01
-7.0604078438e-01
-7.0264686081e-01
-7.0264643835e-01
-7.0261345581e-01
-7.0261291613e-01
-7.0224304298e-01
-7.0223533744e-01
-7.0207505750e-01
-7.0206486461e-01
-7.0104625211e-01
-7.0104594800e-01
-7.0104477640e-01
-7.0104443547e-01
-7.0073999919e-01
-7.0071168445e-01
-7.0069933877e-01
-7.0052638196e-01
-7.0049213802e-01
-7.0047252829e-01
-7.0047134702e-01
-7.0045188820e-01
-7.0045051255e-01
-6.9810224961e-01
-6.9810220056e-01
-6.9809885743e-01
-6.9809880427e-01
-6.9739914275e-01
-6.9739889416e-01
-6.9739294060e-01
-6.9739269352e-01
-6.9617295650e-01
-6.9617286046e-01
-6.9612548692e-01
-6.9612526921e-01
-6.9610045594e-01
-6.9610024421e-01
-6.9589814141e-01
-6.9586738545e-01
-6.9586700311e-01
-6.9585307871e-01
-6.9585293190e-01
-6.9465082494e-01
-6.9299479830e-01
-6.9137949062e-01
-6.8981501472e-01
-6.8757589005e-01
-6.8740813023e-01
-6.8683632101e-01
-6.8626988661e-01
-6.8585154561e-01
-6.8550693040e-01
-6.8472223884e-01
-6.8435444532e-01
-6.8421383158e-01
-6.8398288560e-01
-6.8345490015e-01
-6.8314656599e-01
-6.8304266155e-01
-6.8284721184e-01
-6.8245500892e-01
-6.8222633288e-01
-6.8180896932e-01
-6.8130470734e-01
-6.8082723004e-01
-6.7786710538e-01
-6.7748931668e-01
-6.7708708157e-01
-6.7693870947e-01
-6.7623769451e-01
-6.7591923539e-01
-6.7574091228e-01
-6.7553124962e-01
-6.7534416601e-01
-6.7495506729e-01
-6.7419915037e-01
-6.7292136391e-01
-6.7223787153e-01
-6.7210053224e-01
-6.7148911775e-01
-6.7002748093e-01
-6.6964662221e-01
-6.6933705721e-01
-6.6930861247e-01
-6.6929186166e-01
-6.6896308823e-01
-6.6890583771e-01
-6.6889689625e-01
-6.6884978440e-01
-6.6846284104e-01
-6.6753075514e-01
-6.6748555046e-01
-6.6734717269e-01
-6.6733277784e-01
-6.6548398507e-01
-6.6350693759e-01
-6.6235700674e-01
-6.6235683644e-01
-6.6153562575e-01
-6.6153555652e-01
-6.6153545316e-01
-6.6153537007e-01
-6.5994685870e-01
-6.5994665095e-01
-6.5994649593e-01
-6.5994629913e-01
-6.5872609427e-01
-6.5844967027e-01
-6.5844932069e-01
-6.5842867435e-01
-6.5842832867e-01
-6.5736658622e-01
-6.5733925632e-01
-6.5733011098e-01
-6.5730291325e-01
-6.5380253080e-01
-6.5379900621e-01
-6.5363143619e-01
-6.5362778765e-01
-6.5192274621e-01
-6.5191799005e-01
-6.5190489906e-01
-6.5190471375e-01
-6.5108852255e-01
-6.5108715867e-01
-6.5108660823e-01
-6.5108545581e-01
-6.5079979305e-01
-6.5079911864e-01
-6.5078640640e-01
-6.5078494062e-01
-6.4876506188e-01
-6.4875005880e-01
-6.4381485627e-01
-6.4359393454e-01
-6.4352048518e-01
-6.4334192386e-01
-6.4141879811e-01
-6.4094349144e-01
-6.4078904329e-01
-6.3949728398e-01
-6.3913995098e-01
-6.3820625787e-01
-6.3807167905e-01
-6.3734728305e-01
-6.3731522201e-01
-6.3682237648e-01
-6.3655943941e-01
-6.3611039331e-01
-6.3550070507e-01
-6.3541724596e-01
-6.3539062210e-01
-6.3408540955e-01
-6.3128819513e-01
-6.3126267686e-01
-6.3120657908e-01
-6.3112707087e-01
-6.3107816147e-01
-6.3101468387e-01
-6.2718124898e-01
-6.2714072317e-01
-6.2713625448e-01
-6.2713442975e-01
-6.2652334173e-01
-6.2651579474e-01
-6.2651345818e-01
-6.2649214214e-01
-6.2544007448e-01
-6.2542131395e-01
-6.2542052932e-01
-6.2541962746e-01
-6.2395195589e-01
-6.2395187715e-01
-6.2395177649e-01
-6.2394839918e-01
-6.1898709535e-01
-6.1898697300e-01
-6.1588571362e-01
-6.1588170032e-01
-6.1582743357e-01
-6.1582584258e-01
-6.1445113688e-01
-6.1445101553e-01
-6.1444965546e-01
-6.1444953277e-01
-6.1404145689e-01
-6.1404128803e-01
-6.1404107788e-01
-6.1404091339e-01
-6.1253780289e-01
-6.1253775019e-01
-6.1251591275e-01
-6.1251585277e-01
-6.1118396321e-01
-6.1118360307e-01
-6.1116103205e-01
-6.1116078807e-01
-6.0400870988e-01
-5.8947513129e-01
-5.8222195845e-01
-5.8222191653e-01
-5.8222179027e-01
-5.8222178008e-01
-5.7878022900e-01
-5.7878010902e-01
-5.7877991702e-01
-5.7877979990e-01
-5.7801052378e-01
-5.7795538683e-01
};

\legend{
{\tiny Original Shifts},
{\tiny Added Shifts},
{\tiny Original Ritz Values},
{\tiny Added Ritz Values}
};

\draw[thick,dashed] (-7.3883091259e-01,-1.0) -- (-7.3883091259e-01,2.0);
\draw[thick]        (-7.1855153121e-01,-1.0) -- (-7.1855153121e-01,2.0);
\draw[thick]        (-6.9827214983e-01,-1.0) -- (-6.9827214983e-01,2.0);
\draw[thick]        (-6.5771338707e-01,-1.0) -- (-6.5771338707e-01,2.0);
\draw[thick]        (-6.1715462431e-01,-1.0) -- (-6.1715462431e-01,2.0);
\draw[thick]        (-5.7659586155e-01,-1.0) -- (-5.7659586155e-01,2.0);

\end{axis}
\end{tikzpicture}
  \caption{}
  \label{fig:graphene_shift_insert_evals}
\end{subfigure}~
\begin{subfigure}[b]{0.49\textwidth}
  \begin{tikzpicture}
\begin{axis} [
  width=\textwidth,
  ybar,
  ymin=0,ymax=28,
  every tick label/.append style={font=\tiny},
  xlabel={\tiny SCF Iteration},
  ylabel={\tiny Wall time / s},
  xtick style={draw=none}
]

\addplot coordinates { 
  (0,1.37077215e+01) 
  (1,2*1.35321598e+01) 
  (2,1.37598113e+01)
  (3,1.38519704e+01)
  (4,1.37252327e+01)
};

\end{axis}
\end{tikzpicture}
  \caption{}
  \label{fig:graphene_shift_insert_time}
\end{subfigure}
\caption{Example shift insertion for the Graphene test case. The last slice of this test case was determined to have
missing eigenvalues per the slice validation scheme at SCF iteration 1. (a) shows the positions of the inserted shifts and the new slices
and eigenvalues produced by this insertion. (b) shows the computation required to perform the first 5 SCF iterations
with this insertion.}
\label{fig:graphene_shift_insert}
\end{figure}
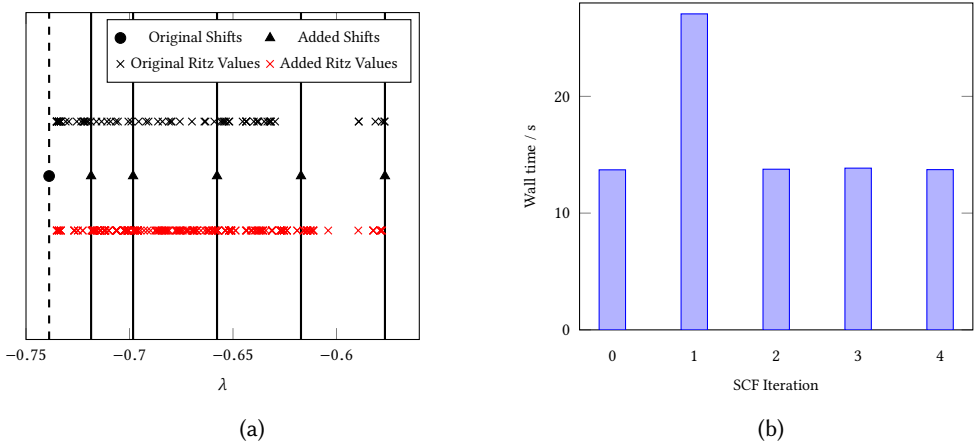

For the Graphene example, we found that some eigenvalues were missed
in the second SCF iteration due to a poorly placed shift produced by
the k-means clustering of approximate eigenvalues obtained in the
first SCF iteration.  The red crosses in \cref{fig:graphene_shift_insert_evals} show
all eigenvalues within the spectral slice $[-0.75, -0.56]$ that contains
the missing eigenvalues.  The black crosses mark the locations of the
approximate eigenvalues that were found by the spectral probe associated
with the poorly placed target shift to the left of this interval (marked
by a solid black circle).

In this case, after performing a DOS estimate of the matrix
pencil updated in the second SCF iteration as suggested
in \cref{sec:misev}, we constructed
five new probes whose target shifts were placed at the positions
marked by the black vertical lines (with arrows) in
\cref{fig:graphene_shift_insert_evals}. After new target shifts were selected from
the k-means clustering of the computed eigenvalues, some of
the inserted probes were mapped to the new shifts and some of
them were deleted in the third SCF iterations. In all subsequent
SCF iterations, no missing eigenvalue was detected, and SCF convergence
was achieved in 10 iterations.


While the inserted new probes captured all the missing eigenvalues, there is an associated cost/penalty for this insertion as shown in \cref{fig:graphene_shift_insert_time}. Because the insertion of new probes essentially
amounts to a recalculation of part of the spectrum in the
second SCF iteration, the wall clock time required to complete that
iteration was doubled.
However, we should point out that this type of probe insertion is rare
in our experiments.  Because it only occurs in early SCF iterations,
the extra cost is typically amortized over the remainder of the SCF
procedure.

\subsection{Probe Basis Dimension}
\label{sec:probe_dim_results}

\begin{figure}[t]

\centering
\begin{subfigure}{0.49\textwidth}
  \centering
  \begin{tikzpicture}
\pgfmathsetlengthmacro\MajorTickLength{
  \pgfkeysvalueof{/pgfplots/major tick length} * 0.5
}
\begin{axis}[
  width=\textwidth,
  legend style={font=\tiny},
  legend cell align={left},
  every tick label/.append style={font=\tiny},
  major tick length=\MajorTickLength,
  ylabel={Wall Time / s},
  xlabel={$k$},
  ymin=0, ymax=5,
  xmin=25,xmax=200,
  x label style={ font=\tiny },
  y label style={ font=\tiny }
]

\addplot[
  mark=o, 
  fill opacity=0,
] table [
  x index=0,
  y expr={\thisrowno{1}/5}
]{
  25  1.57152274e+01
  50  1.67354870e+01
  75  1.72589178e+01
  100 1.76650409e+01
  125 1.85384964e+01
  150 1.93135688e+01
  175 2.07094766e+01
  200 2.14940569e+01
};

\end{axis}
\end{tikzpicture}
  \caption{}
  \label{fig:basis_dim_time}
\end{subfigure}~
\begin{subfigure}{0.49\textwidth}
  \centering
  \begin{tikzpicture}
\pgfmathsetlengthmacro\MajorTickLength{
  \pgfkeysvalueof{/pgfplots/major tick length} * 0.5
}
\begin{axis}[
  width=\textwidth,
  legend style={font=\tiny},
  legend cell align={left},
  every tick label/.append style={font=\tiny},
  major tick length=\MajorTickLength,
  xmin=4, xmax=40,
  ylabel={$\log_{10}(\mathrm{max}\,\|r\|)$},
  xlabel={Subspace Iteration},
  x label style={ font=\tiny },
  y label style={ font=\tiny },
  legend columns=2,
]
\addplot[dashed,mark=triangle*,mark options={solid}, fill opacity=1] table [x=x, y=y1, col sep=comma] {figs/BasisDimConv.dat};
\addplot[mark=o,        fill opacity=0] table [x=x, y=y2, col sep=comma] {figs/BasisDimConv.dat};
\addplot[mark=square*,  fill opacity=1] table [x=x, y=y3, col sep=comma] {figs/BasisDimConv.dat};
\addplot[dashed,mark=square*,mark options={solid},    fill opacity=1] table [x=x, y=y4, col sep=comma] {figs/BasisDimConv.dat};
\addplot[mark=x,        fill opacity=1] table [x=x, y=y5, col sep=comma] {figs/BasisDimConv.dat};
\addplot[mark=triangle, fill opacity=0] table [x=x, y=y6, col sep=comma] {figs/BasisDimConv.dat};
\addplot[mark=x, fill opacity=1] table [x=x, y=y7, col sep=comma] {figs/BasisDimConv.dat};
\addplot[mark=o,dashed, mark options={solid}, fill opacity=1] table [x=x, y=y8, col sep=comma] {figs/BasisDimConv.dat};
\legend{
$k = 25$,
$k = 50$,
$k = 75$,
$k = 100$,
$k = 125$,
$k = 150$,
$k = 175$,
$k = 200$,
};
\end{axis}
\end{tikzpicture}
  \caption{}
  \label{fig:basis_dim_conv}
\end{subfigure}
\caption{
  The effects of the probe basis dimension on timings and convergence in the
  SISLICE method.  Results were obtained using the Graphene test case with 100
  spectral shifts.  Timings (a) were obtained with Intel(R) MKL on 32 Intel(R)
  Haswell threads. The convergence (b) is tracked as the largest residual norm
  for the first spectral slice as a function of subspace iteration.
}
\label{fig:basis_dim}
\end{figure}
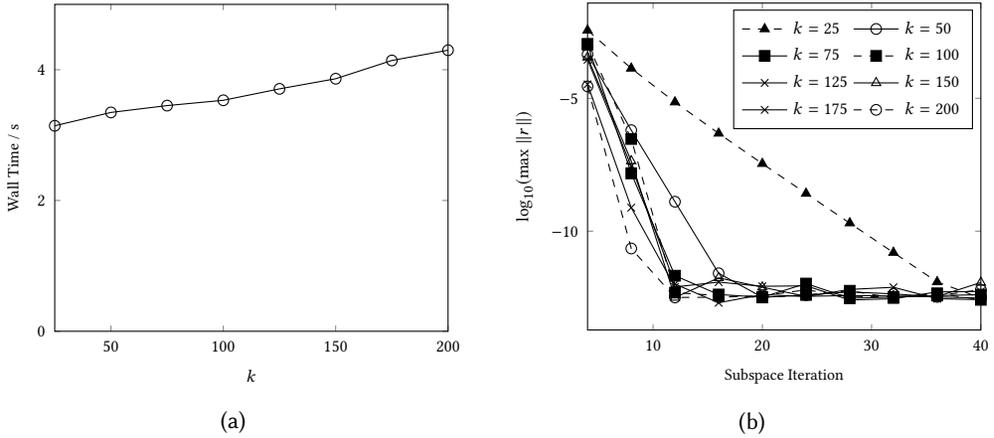

As was discussed in \cref{sec:sisub}, the dimension of the basis used for the
shift-invert subspace iterations need only be \emph{at least} the number of
desired eigenpairs in the neighborhood of a particular spectral shifts. In
practice, the basis dimension should be chosen to be slightly larger to
enable faster convergence.  However, as the basis dimension increases, the
computational time required to perform the shift-invert subspace iterations
also increases due to the need to solve linear systems with a larger number of
right hand sides. In \cref{fig:basis_dim} we examine the effects of basis
dimension on the convergence of the subspace iterations as well as on the
computational time required to perform the subspace iterations for the
Graphene test case.

\Cref{fig:basis_dim_time} shows the increase in computational time required
to perform the subspace iterations as a function of basis dimension.
All timing results were obtained using the Haswell partition of
the Cori Supercomputer (2x16 Intel(R) Xeon(TM) Processor E5-2698 v3 at 2.3
GHz) using Intel(R) MKL to solve the linear systems and are representative of a single set
of 4 subspace iterations. It is clear that even with a 10 fold increase in the basis
dimension, the effect on overall timing for the subspace iterations is negligible.

\Cref{fig:basis_dim_conv} tracks the convergence of a particular spectral slice as
a function of the number of subspace iterations. Due to the nearly uniform
distribution of the eigenvalues within the Graphene's spectrum, this slice is representative of the entire spectrum. The
shift placement was chosen such that each spectral probe is responsible for $\sim 10$
validated eigenpairs. Unlike the effects on timing, the convergence of the examined spectral slice
is sensitive to the basis dimension; with small basis dimensions (e.g. 25) yielding
very sub-optimal convergence results and large basis dimensions (e.g. 100-200) yielding
much faster convergence. In practice, we have found that choosing $k\approx 10\neig/\nsft$
yields sufficiently fast convergence in most cases.

\subsection{Parallel Scalability}
\label{sec:scaling}
In this section, we examine the parallel scaling behavior of the proposed
SISLICE method and detail a performance comparison with the PFEAST method
of \cite{pfeast}.  All timing results were obtained using the Knight's Landing
(KNL) partition of the Cori Supercomputer (Intel(R) Xeon Phi(TM) CPU 7250 @
1.40GHz) and all shifted linear systems (factorizations and backsubstitutions)
were performed using the sparse symmetric solver PARDISO as distributed with
Intel(R) MKL. All shifted inversions were performed on a single KNL node using
64 threads, and all MPI communication is internode.
Although both SISLICE and (P)FEAST
use the shift invert subspace iteration, we will refer to the application of
the contour spectral projector to the basis in the FEAST method as a contour
subspace iteration do differentiate the methods. We note for clarity that this
includes both the complex-shifted linear system solves as well as the
integration reduction required by the FEAST method, i.e. each contour subspace
iteration requires communication of the subspace to some subset of the 
MPI ranks.

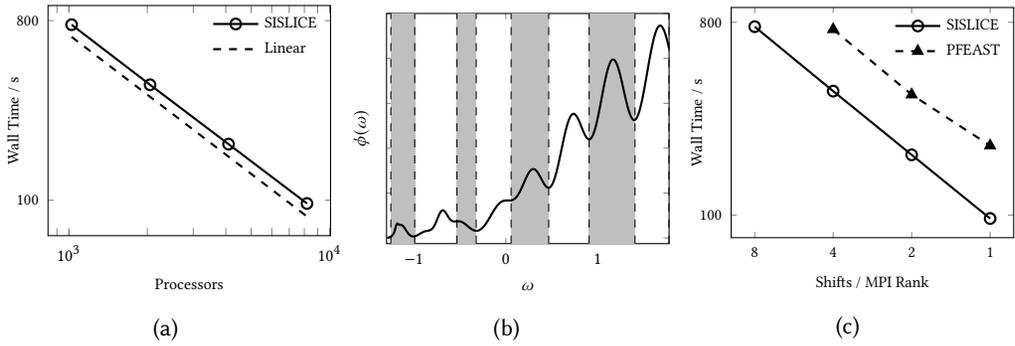
\begin{figure}
\centering
\begin{subfigure}{0.32\textwidth}
  \centering
  \begin{tikzpicture}

\pgfmathsetlengthmacro\MajorTickLength{
  \pgfkeysvalueof{/pgfplots/major tick length} * 0.5
}
\begin{loglogaxis}[
  width=1.2\textwidth,
  every axis plot/.append style={ thick},
  legend style={font=\tiny, draw=none, at={(1,1)}, anchor=north east, fill=none},
  legend cell align={left},
  every tick label/.append style={font=\tiny},
  major tick length=\MajorTickLength,
  xlabel={Processors},
  ylabel={Wall Time / s},
  x label style={ font=\tiny },
  y label style={ font=\tiny },
  ytick={100,800},
  yticklabels={100,800},
  ylabel shift = -8 pt,
]

\addplot[mark=o] table [
  x expr=64*\thisrowno{0},
  y = sislice,
  col sep=comma, 
  unbounded coords=jump] {figs/gaas_sislice_strong_scaling.dat};

\addplot[mark={}, dashed] table [
  x expr=64*\thisrowno{0},
  y expr=(763 - 100) * 16 / \thisrowno{0},
  col sep=comma, 
  unbounded coords=jump] {figs/gaas_sislice_strong_scaling.dat};

\legend{
  SISLICE,
  Linear
};
\end{loglogaxis}
\end{tikzpicture}
  \caption{}
  \label{fig:sislice_strong_scaling}
\end{subfigure}~
\begin{subfigure}{0.32\textwidth}
  \centering
  \begin{tikzpicture}

\pgfmathsetlengthmacro\MajorTickLength{
  \pgfkeysvalueof{/pgfplots/major tick length} * 0.5
}
\begin{axis}[
  width=1.2\textwidth,
  every axis plot/.append style={ thick},
  legend style={font=\tiny},
  legend cell align={left},
  every tick label/.append style={font=\tiny},
  major tick length=\MajorTickLength,
  xmin=-1.3, xmax=1.778,
  ymin=-25, ymax=1000,
  xlabel={$\omega$},
  x label style={ font=\tiny },
  ylabel={$\phi(\omega)$},
  y label style={ font=\tiny },
  ylabel shift = -4 pt,
  ytick={},
  yticklabels={},
]

\addplot +[name path=feast1, dashed, mark={}, black, thin] coordinates {
  (-1.253, \pgfkeysvalueof{/pgfplots/ymin})
  (-1.253, \pgfkeysvalueof{/pgfplots/ymax})
};
\addplot +[name path=feast2, dashed, mark={}, black, thin] coordinates {
  (-9.930e-01, \pgfkeysvalueof{/pgfplots/ymin})
  (-9.930e-01, \pgfkeysvalueof{/pgfplots/ymax})
};
\addplot +[name path=feast3, dashed, mark={}, black, thin] coordinates {
  (-5.330e-01, \pgfkeysvalueof{/pgfplots/ymin})
  (-5.330e-01, \pgfkeysvalueof{/pgfplots/ymax})
};
\addplot +[name path=feast4, dashed, mark={}, black, thin] coordinates {
  (-3.230e-01, \pgfkeysvalueof{/pgfplots/ymin})
  (-3.230e-01, \pgfkeysvalueof{/pgfplots/ymax})
};
\addplot +[name path=feast5, dashed, mark={}, black, thin] coordinates {
  ( 5.700e-02, \pgfkeysvalueof{/pgfplots/ymin})
  ( 5.700e-02, \pgfkeysvalueof{/pgfplots/ymax})
};
\addplot +[name path=feast6, dashed, mark={}, black, thin] coordinates {
  ( 4.670e-01, \pgfkeysvalueof{/pgfplots/ymin})
  ( 4.670e-01, \pgfkeysvalueof{/pgfplots/ymax})
};
\addplot +[name path=feast7, dashed, mark={}, black, thin] coordinates {
  ( 9.070e-01, \pgfkeysvalueof{/pgfplots/ymin})
  ( 9.070e-01, \pgfkeysvalueof{/pgfplots/ymax})
};
\addplot +[name path=feast8, dashed, mark={}, black, thin] coordinates {
  ( 1.407e+00, \pgfkeysvalueof{/pgfplots/ymin})
  ( 1.407e+00, \pgfkeysvalueof{/pgfplots/ymax})
};
\addplot +[name path=feast9, dashed, mark={}, black, thin] coordinates {
  ( 1.777e+00, \pgfkeysvalueof{/pgfplots/ymin})
  ( 1.777e+00, \pgfkeysvalueof{/pgfplots/ymax})
};

\addplot[gray!50] fill between [of=feast1 and feast2];
\addplot[gray!50] fill between [of=feast3 and feast4];
\addplot[gray!50] fill between [of=feast5 and feast6];
\addplot[gray!50] fill between [of=feast7 and feast8];

\addplot[mark={}] table [x=x, y=y, col sep=comma] {figs/gaas_dos.dat};
\end{axis}
\end{tikzpicture}
  \caption{}
  \label{fig:feast_slices}
\end{subfigure}~
\begin{subfigure}{0.32\textwidth}
  \centering
  \begin{tikzpicture}

\pgfmathsetlengthmacro\MajorTickLength{
  \pgfkeysvalueof{/pgfplots/major tick length} * 0.5
}
\begin{loglogaxis}[
  width=1.2\textwidth,
  every axis plot/.append style={ thick},
  legend style={font=\tiny, draw=none, at={(1,1)}, anchor=north east, fill=none},
  legend cell align={left},
  every tick label/.append style={font=\tiny},
  major tick length=\MajorTickLength,
  xlabel={Shifts / MPI Rank},
  ylabel={Wall Time / s},
  x label style={ font=\tiny },
  y label style={ font=\tiny },
  x dir=reverse,
  xtick={1,2,4,8},
  xticklabels={1,2,4,8},
  ytick={100,800},
  yticklabels={100,800},
  ylabel shift = -8 pt,
]

\addplot[mark=o] table [
  x expr=128/\thisrowno{0},
  y=sislice,
  col sep=comma, 
  unbounded coords=jump] {figs/gaas_sislice_strong_scaling.dat};

\addplot[mark=triangle*,dashed,mark options={solid}] table [
  x expr=64/\thisrowno{0},
  y=pfeast,
  col sep=comma, 
  unbounded coords=jump] {figs/gaas_sislice_strong_scaling.dat};

\legend{
  SISLICE,
  PFEAST
};
\end{loglogaxis}
\end{tikzpicture}
  \caption{}
  \label{fig:feast_sislice_compare}
\end{subfigure}
\caption{Scaling of the SISLICE and PFEAST methods. All SISLICE calculations
  were performed to obtain the lowest 1280 eigenvalues of the Ga10As10H30.  (a)
  Strong scaling of SISLICE in comparison with linear scaling. (b) Slice
  intervals for PFEAST.  (c) Strong scaling comparison of SISLICE with PFEAST
  with respect to the number of shifted inversions performed on a single MPI
  rank.}
\label{fig:parallel_scaling}
\end{figure}

\Cref{fig:parallel_scaling} depicts several scaling comparisons between SISLICE
and PFEAST for the Ga10As10 matrix from the SuiteSparse collection
\cite{davis11the} ($N=$ 113,081, $NNZ=$ 6,115,633). Fill reducing reordering
for this matrix was performed using the METIS software \cite{METIS} as packaged
with PARDISO implementation in Intel(R) MKL. Both methods were used to obtain
the lowest 1280 eigenvalues of Ga10As10H30. 
The SISLICE calculations were performed using 128 shifts with $k=100$ and 4
shift-invert subspace iterations per shift. The shifts were selected according
to the DOS-based refinement strategy in \cref{sec:dos}. The PFEAST calculations
were performed with 8 slices using 8 half-contour quadrature points per slice
(a total of 64 total shifted factorizations) and $k=400$. The default elliptic 
contour for used for each slice. As PFEAST does not
admit its own scheme to partition the spectrum of interest into slices, we used
the same DOS partitioning scheme as SISLICE to obtain the intervals for the 8
slices (see \cref{fig:feast_slices}). 8 slices were chosen as to make a fair
comparison with SISLICE by limiting the contours to reasonably localized
eigenvalue clusters as to avoid excessive contour subspace iterations for
embedded spectral regions. PFEAST was run to a trace convergence of $10^{-5}$ to
ensure all eigenpairs were accounted for. All slices converged within 3 contour
subspace iterations.  We note for clarity that the 8 half-contour quadrature
points and $k=400$ basis dimension used per slice is not required for
each of the slices, especially those at the lower end of the spectrum. However,
as the scaling limit of PFEAST will be limited to the computation time required by the
slowest executing contour integral, this does not impact the validity of this
comparison. The basis dimension of 400 was chosen as to allow for convergence
of the larger eigenvalue clusters at the upper end of the considered spectral
domain. However FEAST offers a fallback mechanism to limit the number of
vectors to which the contour subspace iteration is applied for contours
containing only a few eigenvalues. This fallback was triggered for the first
two slices in the FEAST calculations, which produced 21 and 69 vectors after
the first contour subspace iteration, respectively.

The timing given in
\cref{fig:sislice_strong_scaling,fig:feast_sislice_compare} are indicative
all all operations that must take place within a single SCF iteration for the
two methods. This includes the triangular factorizations, backsolves, and
synchronization / integration (as required by FEAST). SISLICE additionally
includes the time required to gather the non-uniformly distributed set of
eigenvectors to a single MPI rank (implemented with \texttt{MPI\_Gatherv}).
Remark that both SISLICE and PFEAST scale linearly out to their strong scaling
limit (1 shift / MPI rank). This is due to the fact that their communication
requirements are minimal and thus pose little to no overhead in execution time.
However, due to the fact that PFEAST uses complex arithmetic, its overall
computational cost per shift is $\sim$2x more than SISLICE. We note for
clarity that due to the differing number of shifts used for SISLICE and PFEAST,
the number of nodes which are represented for the data points in
\cref{fig:feast_sislice_compare} differ by a factor of 2.

\section{Conclusion}
\label{sec:conclusion}

In this work, we have developed the SISLICE method: a robust and efficient
parallel
shift-invert spectrum slicing strategy for self-consistent symmetric eigenvalue
computation. The novelty of the SISLICE method is in its shift selection and
migration strategies which allow for only minimal communication requirements in
its distributed-memory parallel implementation. Like all spectrum slicing methods, the SISLICE
method partitions a spectral region of interest into intervals
which are then treated independently. However, unlike previous and
contemporary slicing methods which rely on effectively sequential shift
placement to partition the spectral region of interest, the SISLICE method
utilizes DOS estimates to form the entire set of spectral slices at once. This
strategy allows for maximal concurrency with minimal communication overhead.
As the desired eigenvalues of the considered matrix pencils
are dynamic throughout the SCF procedure, the SISLICE method employs a shift
migration strategy based on k-means clustering which allows for tracking of the
relevant eigenvalues throughout the SCF procedure without the need to recompute
the costly DOS estimation at each iteration.

We have demonstrated the robustness and parallel efficiency of the SISLICE
method for a representative set of SCF eigenvalue problems commonly encountered
in electronic structure theory in \cref{sec:numerical}.
In particular, we have demonstrated that the k-means shift migration yielded noticeable
convergence improvements in spectral regions with a highly irregular distribution
of eigenvalues (such as the one exhibited for the Silane test case). Further,
even in the cases where it did not yield discernible improvements (e.g. Graphene), its was shown
that the k-means migration did not yield convergence degradation either.

From the perspective of performance, the SISLICE method was demonstrated to
exhibit linear strong scaling for medium to large problem dimensions 
to large processor counts.  Further, we have demonstrated that the main
communication requirement, the synchronization of Ritz values and residual
norms across the distributed network, does not yield a noticeable change in
overall scaling. We have also compared the performance of SISLICE to the 
PFEAST method due to their algorithmic similarity. Both methods scale linearly
out to the strong scaling limit of 1 shift / MPI rank, but the use of complex
arithmetic in PFEAST leads to an overall performance gain of 2x by the SISLICE method.
Further, due to the non-trivial optimal definitions of contours for arbitrary
spectral regions, SISLICE may be treated as a ``black-box" whereas PFEAST requires
considerable tuning to achieve optimal results.

A number of important aspects of the SISLICE method for the SCF
eigenvalue problem, such as the DOS spectral partitioning scheme, k-means
migration strategies and eigenvector seeding for subspace methods, could be
extended to spectrum slicing methods methods such as (P)FEAST, SIPs, and
polynomial filtering in cases where their constituent eigenvalue methods would
be better performing that the shift-invert subspace iteration.  This is of
particular interest for spectra which exhibit similar characteristics as the
all-electron Silane test case which admits several isolated eigenvalue clusters
in the lower region of the spectrum.  Due to the isolated nature of these
clusters, they would likely be better addressed by polynomial filtering,
whereas the larger ``clusters" of eigenvalues higher in the spectrum are well
addressed by SISLICE. This topic will be addressed in future work.

Despite the demonstrated success of the proposed SISLICE method, there are
several topics which were not addressed in this work that should be addressed
to fully demonstrate the effectiveness of the method in real applications.
The first is the integration of the SISLICE method into an actual SCF code such
as NWChemEx, SIESTA, etc. While our results have demonstrated the usefulness
of the SISLICE method for pregenerated matrix sequences, the accuracy of the
eigenvectors at any particular SCF iteration will influence the overall
convergence of the SCF. This topic will also be addressed in future work.


\begin{acks}
This material is based upon work supported by the U.S. Department of Energy,
Office of Science, Office of Advanced Scientific Computing Research, Scientific
Discovery through Advanced Computing (SciDAC) program. (D. B. W. and C. Y.),
and by the Berkeley Lab Undergraduate Research (BLUR) Program, managed by
Workforce Development \& Education at Berkeley Lab (P. B.).  This research was
also partially supported by the Exascale Computing Project (17-SC-20-SC), a
collaborative effort of the U.S. Department of Energy Office of Science and the
National Nuclear Security Administration.  The authors thank the National
Energy Research Scientific Computing (NERSC) center for making computational
resources available to them. They would also like to thank Wei Hu for providing
the Graphene test matrices generated from the SIESTA software and Ajay Panyala
for helping generate the matrix pencils for the Silane example in the NWChemEx
software, as well as Murat Ke{\c{c}}eli for helpful discussions regarding
the development of the SIPs and SIESTA-SIPs software. 
\end{acks}

\bibliographystyle{ACM-Reference-Format}
\bibliography{sislice}

\end{document}